\documentclass[12 pt, draftcls, onecolumn]{IEEEtran}

\usepackage{amssymb,epsfig,color,cite,amsmath,amsfonts,mathrsfs,algorithm,algorithmic}
\usepackage{epstopdf}
\usepackage{graphicx}
\newtheorem{lemma}{Lemma}[section]
\newtheorem{theorem}{Theorem}[section]
\newtheorem{corollary}{Corollary}[section]
\newtheorem{remark}{Remark}[section]

\newtheorem{assumption}{Assumption}[section]

\def\rrr#1\\{\par
\medskip\hbox{\vbox{\parindent=2em\hsize=6.12in
\hangindent=4em\hangafter=1#1}}}
\allowdisplaybreaks

\begin{document}
\title{\huge
Continuous-Time Distributed Algorithms for Extended Monotropic Optimization Problems
}

\author{Xianlin~Zeng, Peng~Yi, Yiguang~Hong, and Lihua~Xie
\thanks{This work was supported by Beijing Natural Science Foundation (4152057), NSFC (61573344), Program 973 (2014CB845301/2/3), and China Postdoctoral Science Foundation (2015M581190).}

\thanks{X.~Zeng, P.~Yi, and Y.~Hong are with Key Laboratory of Systems and Control, Institute of Systems Science, Chinese Academy of Sciences, Beijing, 100190, China (e-mail: xianlin.zeng@amss.ac.cn; yipeng@amss.ac.cn; yghong@iss.ac.cn).}
\thanks{L. Xie is with the School of Electrical and Electronic Engineering, Nanyang Technological University, 50 Nanyang Avenue, Singapore 639798 (e-mail: ELHXIE@ntu.edu.sg).}
}

\maketitle %\baselineskip 14pt \vfill
%\begin{center}
%\today
%\end{center}
%\vspace*{3 em}

%%%%%%%%%%%%%%%%%%%%%%%%%%%%%%%%%%%%%%%%%%%%%%%%%%%%%%%%%%%%%%%%%%%%
%Creating title

%%%%%%%%%%%%%%%%%%%%%%%%%%%%%%%%%%%%%%%%%%%%%%%%%%%%%%%%%%%%%%%%%%%%

%%%%%%%%%%%%%%%%%%%%%%%%%%%%%%%%%%%%%%%%%%%%%%%%%%%%%%%%%%%%%%%%%%%%
%Eliminating page numbering
%\pagestyle{arabic}
%\thispagestyle{arabic}

%\pagenumbering{arabic}
%%%%%%%%%%%%%%%%%%%%%%%%%%%%%%%%%%%%%%%%%%%%%%%%%%%%%%%%%%%%%%%%%%%%

%%%%%%%%%%%%%%%%%%%%%%%%%%%%%%%%%%%%%%%%%%%%%%%%%%%%%%%%%%%%%%%%%%%%

%%%%%%%%%%%%%%%%%%%%%%%%%%%%%%%%%%%%%%%%%%%%%%%%%%%%%%%%%%%%%%%%%%%%
\begin{abstract}\vspace*{.3em}
This paper studies  distributed algorithms for the extended monotropic
optimization problem, which is a general convex optimization problem
with a certain separable structure. The considered objective
function is the sum of local convex functions assigned to agents in
a multi-agent network, with private set constraints and affine
equality constraints. Each agent only knows its local objective function,
local constraint set, and neighbor information. %, in addition to the decomposed information of the problem.
We propose two novel
continuous-time distributed subgradient-based algorithms with
projected output feedback and  derivative feedback,
respectively, to solve the extended monotropic optimization problem.
Moreover, we show that the algorithms converge to the optimal
solutions under some mild conditions, by virtue of variational
inequalities, Lagrangian methods, decomposition methods, and
nonsmooth Lyapunov analysis.   Finally, we give two examples to
illustrate the applications of the proposed algorithms.

%distributed method for convex optimization problems with a certain separability structure
\vspace{3em}


\vspace{1em}

\end{abstract}

\IEEEpeerreviewmaketitle

\baselineskip=22pt

%\clearpage
\section{Introduction}

Recently, distributed (convex) optimization problems have been
studied in many fields including sensor networks, neural learning
systems, and power systems
\cite{LW:TAC:2015,BX:2009,YHL:SCL:2015,YHX:TNNLS:2016}.
Such problems are often formulated in terms of an objective function
in the form of a sum of individual objective functions, each of which
represents the contribution of an agent to the global objective function.
Various distributed approaches for efficiently solving convex
optimization problems in multi-agent networks have been proposed
since these distributed algorithms have advantages over centralized
ones for  large-scale optimization problems, with inexpensive and
low-performance computations for each agent/node.  Both
continuous-time and discrete-time algorithms for distributed
optimization without any constraint were extensively studied
recently \cite{NOP:2010,GC:2014,SJH:TAC:2013}.

It is known that continuous-time algorithms are becoming more and more
popular, especially when the optimization may be achieved by
continuous-time physical systems. In fact, the neurodynamic approach
is one of the continuous-time optimization approaches, and has been well developed for numerous neural network models
\cite{LYW:NNLS:2016,TH:1986,CHZ:1999,LW:MAY2013}; and the optimization in smart grids has also been studied based on continuous-time
dynamical systems \cite{GC:2014,YHL:SCL:2015}.

In practice, many optimization problems get involved with various
constraints and nonsmooth objective functions
\cite{PSL:2003,LW:TAC:2015,YHL:SCL:2015,GC:2014}.
Subgradients and projections have been found to be widely used in the
study of distributed nonsmooth optimization with constraints
\cite{YHL:SCL:2015,RM:2015,arXiv:YHL}. In fact, different algorithms were proposed in the form of differential inclusions to solve nonsmooth optimization problems in the literature (see \cite{BX:2009,FNQ:2004,LW:MAY2013}).

The monotropic optimization is an optimization problem with
separable objective functions, affine equality constraints, and
set constraints.  As a generalization of linear programming and
network programming, it was first introduced and extensively studied
by \cite{Rockafellar:1984,Rockafellar1985}.   The extended
monotropic optimization (EMO) problem was then studied in
\cite{Bertsekas:JOTA:2008}, but not in a distributed manner.
%{\color{red}which is?? related to?? (please do not write several
%``applications"s, which are not the same, in one short paragraph??)
%the resource allocation problems
%\cite{XJB:2004,YHL:SCL:2015,arXiv:YHL}, network utility problems
%\cite{SS:2008}, and network flow problems
%\cite{Rockafellar:1984,Rockafellar1985}.}
In fact, in many applications
such as wireless communication, sensor networks, neural computation, and  networked
robotics, the optimization problems can be converted to EMO problems
\cite{Bertsekas:JOTA:2008}.  Recently, distributed algorithms
and decomposition methods for EMO problems have become more and more
important, as pointed out in \cite{CDZ:MPA:2015}.  However, due to the
complicatedness of EMO problems, very few distributed optimization designs for them  have been proposed.

This paper aims  to investigate EMO problems  with
nonsmooth objective functions via distributed approaches.
Based on Lagrangian functions and decomposition methods, our
distributed algorithms are given based on projections to achieve the
optimal solutions of the problems. The analysis of the proposed
algorithms is carried out by applying the stability theory of differential
inclusions to tackle nonsmooth objective functions. The main
technical contributions of the paper are three folds. Firstly, we
study the distributed EMO problem, and propose two novel distributed
continuous-time algorithms for the problem, by using projected
output feedback and derivative feedback, respectively, to
deal with the nonsmoothness and set constraints.  Secondly, based on
the Lagrangian function method along with projections, we show that
the proposed algorithms have bounded states and  solve the  EMO
problems with any initial condition, by virtue of the
stability theory of differential inclusions and nonsmooth analysis
techniques. Thirdly, we apply the given algorithms to two practical
problems and illustrate their effectiveness.

The rest of the paper is organized as follows:  Section
\ref{sec:def} shows the preliminary knowledge related to graph
theory, nonsmooth analysis, convex optimization, and projection
operators.  Next, Section \ref{Distributed_Optimization} formulates
a class of distributed extended monotropic optimization (EMO) problems with
nonsmooth objective functions, while Section \ref{algorithm}
proposes two distributed algorithms based on projected output
feedback and  derivative feedback.  Then Section \ref{main_result}
shows the convergence of the proposed algorithms with nonsmooth
analysis. Following that, Section \ref{sec_appl} gives two
application examples to illustrate the theoretical results. Finally,
Section \ref{conclusion} presents some concluding remarks.

\section{Mathematical Preliminaries}

In this section, we introduce relevant notations, concepts, and
preliminaries on graph theory, differential inclusions,
convex analysis, and projection operators.

\label{sec:def}
\subsection{Notation}
$\mathbb R$ denotes the set of real numbers;
$\mathbb R^n$ denotes the set of $n$-dimensional real column
vectors; $\mathbb R^{n\times m}$ denotes the set of $n$-by-$m$ real
matrices; $I_n$ denotes the $n\times n$ identity matrix; and
$(\cdot)^\mathrm{T}$ denotes the transpose. Denote $\mathrm
{rank}\,A$ as the rank of the matrix $A$, $\mathrm{range} (A)$ as
the range of $A$,
$\ker (A)$ as the kernel of $A$, $\mathrm{diag}\{A_1,\ldots,A_n\}$ as the block diagonal matrix of $A_1,\ldots,A_n$, %$\lambda _{\mathrm {max}} (A)$ for the largest eigenvalue of the matrix $A$,
$ 1_{n}$ ($1_{n\times q}$) as the $n\times 1$ vector ($n\times q$
matrix) with all elements of 1, $ 0_{n}$ ($0_{n\times q}$) as the
$n\times 1$ vector ($n\times q$ matrix) with all elements of 0, and
$A \otimes B$ as the Kronecker product of matrices $A$ and $B$. $A>
0$ ($A\geq 0$) means that matrix $A\in\mathbb R^{n\times n}$ is
positive definite (positive semi-definite).  Furthermore,
$\|\cdot\|$ stands for the Euclidean norm;
%$\|\cdot\|_p$ denotes the $p$-norm where $p\geq 1$,
$\overline{\mathcal{S}}$ ($\mathcal S^\circ$) for the closure
(interior) of the subset $\mathcal{S}\subset\mathbb{R}^{n}$;
%$\mathrm{int}(\mathcal{S})$ denotes the interior of the subset $\mathcal{S}$, $\mathrm{dim}(\mathcal S)$ denotes the dimension of the vector space $\mathcal S$,
$\mathcal B_{\epsilon}(x),x\in\mathbb R^n,\epsilon>0,$ for the open
ball {\em centered} at $x$ with {\em radius} $\epsilon$.
$\mathrm{dist}(p,\mathcal M)$ denotes the distance from a point $p$
to the set $\mathcal M$ (that is, $\mathrm{dist}(p,\mathcal M)
\triangleq \inf_{x\in\mathcal M}\|p-x\|$), and $x(t)\rightarrow
\mathcal M$ as $t\rightarrow \infty$ denotes that $x(t)$ approaches
the set $\mathcal M$ (that is, for each $\epsilon>0$, there is $T>0$
such that $\mathrm{dist}(x(t),\mathcal M)<\epsilon$ for all $t>T$).

\subsection{Graph Theory}
A weighted undirected graph is described by $\mathcal G$ or
$\mathcal G(\mathcal V, \mathcal E, A)$, where $\mathcal
V=\{1,\ldots, n\}$ is the set of nodes, $\mathcal E\subset\mathcal V
\times \mathcal V$ is the set of edges, $ A=[a_{i,j}]\in\mathbb
R^{n\times n}$ is the {\em weighted adjacency matrix} such that
$a_{i,j}=a_{j,i}>0$ if $\{j,i\}\in\mathcal E$ and $a_{i,j}=0$
otherwise. The {\em weighted Laplacian matrix} is $L_n=D- A$, where
$D\in\mathbb R^{n\times n}$ is diagonal with $D_{i,i}=\sum_{j=1}^n
a_{i,j}$, $i\in\{1,\ldots,n\}$. In this paper, we call $L_n$ the
Laplacian matrix and $A$ the adjacency matrix of $\mathcal G$ for
convenience when there is no confusion. Specifically, if the
weighted undirected graph $\mathcal G$ is connected, then
$L_n=L_n^{\rm T}\geq 0$, $\mathrm{rank} \,L_n=n-1$ and $\ker
(L_n)=\{k{1}_n:k\in\mathbb R\}$ \cite{GR2001}.

\subsection{Differential Inclusion}

Following \cite{AC:1984}, a differential inclusion is
given by
\begin{eqnarray}\label{DI}
  \dot x(t) \in  \mathcal F (x(t)),\quad x(0)=x_0,\quad t\geq 0,
\end{eqnarray}
where $\mathcal F$ is a set-valued map from $\mathbb R^q$ to the
compact, convex subsets of $\mathbb R^q$. For each state
$x\in\mathbb{R}^{q}$, system \eqref{DI} specifies a
\textit{set} of possible evolutions  rather than a
single one.
A %Caratheodory
solution of \eqref{DI} defined on $[0, \tau] \subset [0,\infty)$ is
an absolutely continuous function $x:[0, \tau]\rightarrow \mathbb
R^q$ such that \eqref{DI} holds for almost all $t\in [0, \tau]$ for
$\tau>0$. The solution $t\mapsto x(t)$ to (\ref{DI}) is a
\textit{right maximal solution} if it cannot be extended forward in
time. Suppose that all right maximal solutions to (\ref{DI}) exist
on $[0,\infty)$.
A set $\mathcal M$ is said to be {\em weakly invariant}
(resp., {\em strongly invariant}) with respect to (\ref{DI}) if, for every $x_0\in\mathcal M$, $\mathcal M$ contains a maximal solution %\cite{R:1998,BC:1999}
(resp., all maximal solutions) of (\ref{DI}).
A point $z\in\mathbb R^{q}$ is a {\em positive limit point} of a solution $\phi(t)$ to (\ref{DI}) with $ \phi(0) = x_0\in\mathbb R^q$, if there
exists a sequence $\{t_k\}_{k=1}^{\infty}$ with $t_k\rightarrow \infty$ and $\phi(t_k)\rightarrow z$ as $k\rightarrow\infty$. The set $\omega(\phi(\cdot))$ of all such positive limit points is the {\em positive limit set} for the trajectory $\phi(t)$ with $ \phi(0) = x_0\in\mathbb R^q$.

An equilibrium point of (\ref{DI}) is a point $x_e\in\mathbb R^q$ such that ${0}_q\in\mathcal F(x_e)$. It is easy to see that $x_e$ is an equilibrium point of (\ref{DI}) if and only if the constant function $x(\cdot)=x_e$ is a solution of (\ref{DI}). An equilibrium point $z\in\mathbb R^q$ of (1) is {\em Lyapunov stable} if, for every $\varepsilon>0$, there exists $\delta=\delta(\varepsilon)>0$ such that, for every initial condition $x(0)=x_0\in\mathcal B_{\delta}(z)$, every solution $x(t)\in\mathcal B_{\varepsilon}(z)$ for all $t\geq 0$.

Let $V:\mathbb R^q\rightarrow \mathbb R$ be a locally Lipschitz
continuous function, and $\partial{V}$  the \textit{Clarke
generalized gradient} \cite{Clarke:1983} of $V(x)$ at $x$. The
\textit{set-valued Lie derivative} \cite{Clarke:1983}
$\mathcal{L}_{\mathcal
F}V:\mathbb{R}^{q}\to\mathfrak{B}(\mathbb{R})$ of $V$ with respect
to (\ref{DI}) is defined as $\mathcal{L}_{\mathcal F}
V(x)\triangleq\{a\in\mathbb{R}:{\rm{there}}\,\,{\rm{exists}}\,\,v\in\mathcal{F}(x)\,\,{\rm{such}}\,\,{\rm{that}}\,\,p^{\rm{T}}
v=a\,\,{\rm{for}}\,\,{\rm{all}}\,\,p\in\partial{V}(x)\}.$ In the
case when $\mathcal{L}_{\mathcal F}V(x)$ is nonempty, we use
$\max\mathcal{L}_{\mathcal F}V(x)$ to denote the largest element of
$\mathcal{L}_{\mathcal F}V(x)$. Recall from reference \cite{BC:1999}
that, if $\phi(\cdot)$ is a solution of (\ref{DI}) and
$V:\mathbb{R}^{q}\to\mathbb{R}$ is locally Lipschitz and {\em
regular} (see \cite[p.~39]{Clarke:1983}), then $\dot V(\phi(t))$
exists almost everywhere, and $\dot
V(\phi(t))\in\mathcal{L}_{\mathcal F}V(\phi(t))$ almost everywhere.

Next, we introduce a version of the invariance principle (Theorem 2 of \cite{Cortez:2008}), which is based on nonsmooth regular functions.

\begin{lemma}\cite{Cortez:2008}\label{nonsmooth_invariance}
For the differential inclusion (\ref{DI}), we assume that $\mathcal
F$ is upper semicontinuous and locally bounded, and $\mathcal F(x)$
takes nonempty, compact, and convex values. Let
$V:\mathbb{R}^{q}\to\mathbb{R}$ be a locally Lipschitz and regular
function, $\mathcal S\subset\mathbb{R}^q$ be compact and strongly
invariant for (\ref{DI}), $\phi(\cdot)$ be a solution of
(\ref{DI}),
$$
\mathcal R=\{x\in\mathbb R^q:0\in\mathcal{L}_{\mathcal F}V(x)\},
$$
and $\mathcal M$ be the largest weakly invariant subset of
$\overline{\mathcal R}\cap \mathcal S$, where $\overline{\mathcal R}$ is the closure of ${\mathcal R}$. If $\max \mathcal
{L}_{\mathcal F} V(x)\leq 0$ for all $x\in \mathcal S$, then
$\mathrm{dist}(\phi(t),\mathcal M)\rightarrow 0$ as $t\rightarrow
+\infty$.
\end{lemma}

\subsection{Convex Analysis}
A function $\psi:\mathbb R^q\rightarrow \mathbb R$ is {\em convex}
if $\psi(\lambda x+(1-\lambda)y)\leq \lambda
\psi(x)+(1-\lambda)\psi(y)$ for all $x,y\in\mathbb{R}^q$ and
$\lambda\in [0,1]$. A function $\psi:\mathbb R^q\rightarrow \mathbb
R$ is {\em strictly convex} whenever $\psi(\lambda
x+(1-\lambda)y)<\lambda \psi(x)+(1-\lambda)\psi(y)$ for all
$x,y\in\mathbb{R}^q$, $x\not= y$ and $\lambda\in (0,1)$.
%A convex function $f$ is strongly convex if, and only if, there exists a constant $\delta>0$ such that $f(x)-\delta \|x\|^2$ is convex.
%The effective domain of $\psi$ is defined by $\mathrm{dom}\,\psi\triangleq\{x\in\mathbb R^q: \psi(x)<+\infty\}$. A convex function $\psi$ is said to be {\em proper} \cite[p.~530]{ZSN:2003} if $\mathrm{dom}\,\psi\not=\emptyset$ and $\psi>-\infty$ for all $x\in\mathrm{dom}\,\psi$.
%Let $\psi:\mathbb R^q\rightarrow \mathbb R\bigcup\{+\infty\}$ be a proper convex function.
Let $\psi:\mathbb R^q\rightarrow \mathbb R$ be a convex function. The {\em subdifferential} \cite[p.~544]{ZSN:2003} $\partial_{\rm sub} \psi$ of $\psi$ at $x\in\mathbb R^q$ is defined by
$
\partial_{\rm sub}{\psi}(x)\triangleq \{p\in\mathbb R^q: \langle p,y-x\rangle\leq {\psi(y)-\psi(x)} ,\,\forall y\in\mathbb R^q\},
$
and the elements of $\partial_{\rm sub}{\psi}(x)$ are called {\em subgradients} of $\psi$ at point $x$. %Note that continuous convex functions are locally Lipschitz continuous \cite[p.~600]{ZSN:2003}, regular \cite[p.~608]{ZSN:2003}, and their subdifferentials and Clarke generalized gradients coincide \cite[p.~607]{ZSN:2003}.
Recall from  \cite[p.~607]{ZSN:2003}  that continuous convex functions are locally Lipschitz
continuous, regular, and their subdifferentials and Clarke
generalized gradients coincide. Thus, the framework for stability
theory of differential inclusions can be applied to the theoretical
analysis in this paper.

The following result can be easily verified by the property of strictly convex functions.%is immediate by extending Theorem 2.67 of \cite{Ruszczynski:2006}.

\begin{lemma}\label{convex_strict}
Let $f:\mathbb{R}^q\rightarrow\mathbb{R}$ be a continuous strictly
convex function.
%  that is,
%  $$f(\lambda x+(1-\lambda)y)<\lambda f(x)+(1-\lambda)f(y)$$
%  for all $x,y\in\mathbb{R}^q$, $x\not= y$ and $\lambda\in (0,1)$.
  Then
  \begin{eqnarray}\label{strictly_convex}
  % \nonumber % Remove numbering (before each equation)
    (g_x-g_y)^{\rm T}(x-y)>0,
    %(g_x-g_y)^{\rm T}(x-y)>0,\quad x,y\in\mathbb{R}^q, \quad x\not= y,\quad g_x\in\partial f(x),\quad g_y\in\partial f(y).
  \end{eqnarray}
  for all $x\not= y$, where $g_x\in\partial f(x)$ and $g_y\in\partial f(y)$.
\end{lemma}

\subsection{Projection Operator}
Define $P_{\Omega}(\cdot)$ as a projection operator given by $P_{\Omega}(u)=\mathrm{arg}\,\min_{v\in\Omega}\|u-v\|$, where $\Omega\subset \mathbb R^n$ is closed and convex. A basic property \cite{KS:1982}
of a projection $P_{\Omega}(\cdot)$ on a closed convex set $\Omega\subset \mathbb R^n$ is %\label{Projection_lemma}
\begin{eqnarray}\label{projection_inequality}
(u-P_{\Omega}(u))^{\rm T}(v-P_{\Omega}(u))\leq 0,\quad \forall u\in\mathbb R^n,\quad \forall v\in\Omega.
\end{eqnarray}

Using (\ref{projection_inequality}), the following results can be easily verified.
\begin{lemma}\label{lemma_ineq}
If $\Omega\subset \mathbb R^n$ is closed and convex, then $(P_{\Omega}(x)-P_{\Omega}(y))^{\rm T}(x-y)\geq \|P_{\Omega}(x)-P_{\Omega}(y)\|^2$ for all $x,y\in\mathbb{R}^{n}$.
\end{lemma}
%\begin{IEEEproof}
%It follows from (\ref{projection_inequality}) that
%\begin{eqnarray}
%% \nonumber % Remove numbering (before each equation)
%  (x-P_{\Omega}(x))^{\rm T}(P_{\Omega}(y)-P_{\Omega}(x)) &\leq & 0, \label{proj_ineq1}\\
%  (y-P_{\Omega}(y))^{\rm T}(P_{\Omega}(x)-P_{\Omega}(y)) &\leq & 0, \label{proj_ineq2}
%\end{eqnarray}
%where $x,y\in\mathbb{R}^n$.
%
%By adding \eqref{proj_ineq1} and \eqref{proj_ineq2}, it follows that $(P_{\Omega}(x)-P_{\Omega}(y))^{\rm T}(x-y)\geq \|P_{\Omega}(x)-P_{\Omega}(y)\|^2$.
%\end{IEEEproof}

\begin{lemma}\label{projection_diff}\cite{LW:MAY2013}
Let $\Omega\subset \mathbb R^n$ be closed and convex, and define $V:\mathbb{R}^n\rightarrow\mathbb{R}$ as $V(x)=\frac{1}{2}(\|x-P_{\Omega}(y)\|^2-\|x-P_{\Omega}(x)\|^2)$ where $y\in\mathbb{R}^n$. Then $V(x)\geq \frac{1}{2}\|P_{\Omega}(x)-P_{\Omega}(y)\|^2$,  $V(x )$ is differentiable and convex with respect to $x$, and $\nabla V(x)= P_{\Omega}(x)-P_{\Omega}(y)$.
\end{lemma}

\section{Problem Description }
\label{Distributed_Optimization}
In this section, we present the distributed extended monotropic optimization (EMO) problem
with nonmsooth objective functions, and give the optimality condition for the problem.

Consider a network of $n$ agents interacting over a graph $\mathcal G$.
There are a local objective function $f^{i}:\Omega_i\rightarrow\mathbb R$ and a local feasible constraint set $\Omega_i\subset\mathbb R^{q_i}$ for all $i\in\{1,\ldots,n\}$. Let $x_i\in \Omega_i\subset\mathbb R^{q_i}$ and denote $x\triangleq [x_1^{\rm T},\ldots,x_n^{\rm T}]^{\rm T}\in\Omega\triangleq \prod_{i=1}^n \Omega_i\subset\mathbb R^{\sum_{i=1}^{n}q_i}$.
The global objective function of the network is $f(x) = \sum_{i=1}^n f^i(x_i),\,x\in\Omega\subset\mathbb R^{\sum_{i=1}^{n}q_i}$.

Here we consider the following distributed EMO problem
\begin{subequations}\label{optimization_problem}
\begin{align}
&\min\,f(x), \quad f(x)= \sum_{i=1}^n f^i(x_i),\\
 & Wx=\sum_{i=1}^n W_i x_i = d_0,\quad x_i\in \Omega_i\subset\mathbb R^{q_i},\quad i\in\{1,\ldots,n\},
\end{align}
\end{subequations}
where  $W_i\in\mathbb{R}^{m\times q_i},\,i\in\{1,\ldots,n\}$ and $W = [W_1,\ldots,W_n]\in\mathbb{R}^{m\times \sum_{i=1}^{n}q_i}$. In this problem, agent $i$ has its state $x_i\in \Omega_i\subset\mathbb R^{q_i}$, objective function $f_i(x_i)$, set constraint $\Omega_i\subset\mathbb R^{q_i}$, constraint matrix $W_i\in\mathbb{R}^{m\times q_i}$, and information from neighboring agents.

%Consider a network of $n$ agents interacting over a graph $\mathcal G$ and the following distributed EMO problem
%\begin{subequations}\label{optimization_problem}
%\begin{align}
%&\min\,f(x), \quad f(x)= \sum_{i=1}^n f^i(x_i),\\
% & Wx=\sum_{i=1}^n W_i x_i = d_0,\quad x_i\in \Omega_i\subset\mathbb R^{q_i},\quad i\in\{1,\ldots,n\},
%\end{align}
%\end{subequations}
%where agent $i$ has its state $x_i\in \Omega_i\subset\mathbb R^{q_i}$, local objective function $f_i:\Omega_i\rightarrow\mathbb R$, local set constraint $\Omega_i\subset\mathbb R^{q_i}$, constraint matrix $W_i\in\mathbb{R}^{m\times q_i}$, and information from neighboring agents. Denote $x\triangleq [x_1^{\rm T},\ldots,x_n^{\rm T}]^{\rm T}\in\Omega\triangleq \prod_{i=1}^n \Omega_i\subset\mathbb R^{\sum_{i=1}^{n}q_i}$ and $W = [W_1,\ldots,W_n]\in\mathbb{R}^{m\times \sum_{i=1}^{n}q_i}$.
%The global objective function of the network is $f(x) = \sum_{i=1}^n f^i(x_i),\,x\in\Omega\subset\mathbb R^{\sum_{i=1}^{n}q_i}$.

The goal of the distributed  EMO is to solve the problem in a distributed manner. In a distributed optimization algorithm, each agent  in the graph $\mathcal G$ only uses its own local cost function, its local set constraint, the decomposed information of the global equality constraint, and the shared information
of its neighbors through constant local communications.
The special case of problem \eqref{optimization_problem}, where each component $x_i$ is one-dimensional (that is, $q_i = 1$), is called the monotropic programming problem and  has been introduced and studied extensively in \cite{Rockafellar:1984,Rockafellar1985}.

\begin{remark}
The distributed EMO problem (\ref{optimization_problem}) covers many problems in the recent distributed optimization studies because of the general expression. For
example, it generalizes the optimization model in resource allocation problems \cite{YHL:SCL:2015,arXiv:YHL} by allowing nonsmooth objective functions and a more general equality constraint.
Moreover, it covers the model proposed in \cite{LYW:NNLS:2016} and generalizes the model in the distributed constrained optimal consensus problem \cite{QLX:A:2016} by allowing heterogeneous constraints.
\end{remark}

For illustration, we introduce two special cases of our problem:
\begin{itemize}
%\item Consider the following optimization problem \cite{Kia:CDC:2015}, which can be viewed as an extension of the conventional resource allocation problem.
%\begin{eqnarray}\label{special_multiple_constraint}
% \min f(x),\quad  f(x)=\sum_{i=1}^n f^i(x_i), \quad \sum_{i=1}^n \omega^j_i x_i=b^j,\quad j=1,\ldots,m,
%\end{eqnarray}
%where $x\triangleq [x_1^{\rm T},\ldots,x_n^{\rm T}]^{\rm T}\in\mathbb R^{nq}$, $x_i\in\mathbb R^q$, $\omega^j_i\in\mathbb R$, $b^j\in\mathbb R^q$, and $f^i:\mathbb R^q\rightarrow\mathbb R$ for $i\in\{1,\ldots,n\}$ and $j\in\{1,\ldots,m\}$. Each agent $i$ only knows state $x_i$, local objective function $f^i(x_i)$, parameter $\omega^j_i$ and neighboring information for $i\in\{1,\ldots,n\}$ and $j\in\{1,\ldots,m\}$.
%Let $d_0\triangleq [b_1^{\rm T},\ldots,b_m^{\rm T}]^{\rm T}$, let $W_i\triangleq \begin{bmatrix}\omega_i^1 I_q\\
%\vdots\\
%\omega_i^m I_q
%\end{bmatrix},\,i\in\{1,\ldots,n\}$, and define $W \triangleq [W_1,\ldots,W_n]\in\mathbb{R}^{q\times nq}$. Problem \eqref{special_multiple_constraint} can be written in the form of (\ref{optimization_problem}). Hence, EMO problem described by (\ref{optimization_problem}) covers the optimization problem with several equality constraints given by \eqref{special_multiple_constraint}.

\item Consider the following optimization problem investigated in \cite{LYW:NNLS:2016}
\begin{subequations}\label{local_global_constraint}
\begin{align}
% \nonumber % Remove numbering (before each equation)
 & \min  f(x),\quad  f(x)=\sum_{i=1}^n f^i(x_i), \quad x_i\in\Omega_i=\{x_i\in\mathbb R^q:g_i(x_i)\leq 0\},\label{local_global_constraint_1}\\
  & A_i x_i=b_i,\quad i\in\{1,\ldots,n\}, \label{local_global_constraint_2}\\
  & Lx = 0_{nq}, \label{local_global_constraint_3}
\end{align}
\end{subequations}
where $x\triangleq [x_1^{\rm T},\ldots,x_n^{\rm T}]^{\rm T}\in\mathbb R^{nq}$, $A_i\in\mathbb R^{m_i\times q}$, $b_i\in\mathbb R^{m_i}$, $L\in\mathbb R^{nq\times nq}$, $x_i\in\mathbb R^q$,   and $f^i:\mathbb R^q\rightarrow\mathbb R$ for all $i\in\{1,\ldots,n\}$. Equation \eqref{local_global_constraint_2} is the local equality constraint for agent $i$, and equation \eqref{local_global_constraint_3} is the global equality constraint. %Agent $i$ only knows state $x_i$, local objective function $f^i(x_i)$, matrices $A_i$, $B_i$, and the states of its neighboring agents for $i\in\{1,\ldots,n\}$.
Let $d_0\triangleq [b_1^{\rm T},\ldots,b_n^{\rm T}, 0_{nq}^{\rm T}]^{\rm T}$, $A=\mathrm{diag} [A_1,\dots,A_n]\in\mathbb{R}^{\sum_{i=1}^{n}m_i\times nq}$, and $W\triangleq \begin{bmatrix}A \\L\end{bmatrix}\in\mathbb{R}^{(nq+\sum_{i=1}^{n}m_i)\times nq}$.
Problem \eqref{local_global_constraint} can be written in the form of (\ref{optimization_problem}). Hence, the EMO problem described by (\ref{optimization_problem}) covers the optimization problem with local and global equality constraints given by \eqref{local_global_constraint}.

\item Consider the minimal norm problem of underdetermined linear equation \cite{arXiv:0905.4745}
\begin{eqnarray}\label{local_global_constraint2}
% \nonumber % Remove numbering (before each equation)
  \min  f(x),\quad  f(x)=\sum_{i=1}^n   \|x_i \|^2, \quad \sum_{i=1}^n A_i x_i=b,
\end{eqnarray}
where $x\triangleq [x_1^{\rm T},\ldots,x_n^{\rm T}]^{\rm T}\in\mathbb R^{\sum_{i=1}^{n}q_i}$, $A_i\in\mathbb R^{m\times q_i}$, $b\in\mathbb R^{m}$, $x_i\in\mathbb R^{q_i}$, and $f^i:\mathbb R^{q_i}\rightarrow\mathbb R$ for all $i\in\{1,\ldots,n\}$. Hence, the EMO problem also covers the minimal norm problem of linear algebraic equation. In the description of our problem, each agent knows a block $A_i$, which is different from the framework in \cite{MLM:TAC:2015} solving the distributed linear equation, where each agent knows a subset of the rows of $[A_1,\ldots,A_n]$ and $b$. If the number of variables is large, our framework  obviously has advantages on reducing the computation and information  loads of agents.

\end{itemize}

To ensure the well-posedness of the problem, the following assumption for problem (\ref{optimization_problem}) is needed, which is quite standard.

\begin{assumption}\label{Assumption}
\begin{enumerate}
  \item The weighted graph $\mathcal G$ is connected and undirected.
  \item For all $i\in\{1,\dots,n\}$, $f^i$ is %continuous and
  strictly convex on an open set containing $\Omega_i$, and $\Omega_i\subset\mathbb R^{q_i}$ is {closed} and {convex}.
  \item (Slater's constraint condition) There exists $x\in\Omega^\circ$ satisfying the constraint $Wx=d_0$, where $\Omega^\circ$ is the interior of $\Omega$.
\end{enumerate}
\end{assumption}

\begin{lemma}\label{equiv_cond}
Under Assumption \ref{Assumption}, ${x}^*\in\Omega$ is an optimal solution of (\ref{optimization_problem}) if and only if there exist $\lambda^*_0\in\mathbb R^{m}$ and $g(x^*)\in\partial f(x^*)$ such that
\begin{eqnarray}
  && x^*=P_{\Omega}(x^*  -  g({x^*})+W^{\rm T}\lambda_0^*), \label{opt_equi_1}\\
  && Wx^*=d_0.\label{opt_equi_2}
\end{eqnarray}
\end{lemma}

\begin{IEEEproof}
Consider problem (\ref{optimization_problem}).  By the KKT optimality condition (Theorem 3.25 of \cite{Ruszczynski:2006}),  ${x}^*\in\Omega$ is an optimal solution of (\ref{optimization_problem}) if and only if  there exist $\lambda^*_0\in\mathbb R^{m}$ and $g(x^*)\in\partial f(x^*)$ such that \eqref{opt_equi_2} holds and
\begin{eqnarray}\label{equi_1}
-g({x^*})+W^{\rm T}\lambda_0^*\in \mathcal N_{\Omega}(x^*),
\end{eqnarray}
where $\mathcal N_{\Omega}(x^*)$ is the normal cone of $\Omega$ at an element $x^*\in\Omega$. Note that \eqref{equi_1} holds if and only if \eqref{opt_equi_1} holds. Thus, the proof is completed.
\end{IEEEproof}

\section{Optimization Algorithms}
\label{algorithm}
In this section, we propose two distributed optimization algorithms to solve the EMO problem with nonsmooth objective functions. To our best knowledge, there are no distributed continuous-time algorithms for such problems with rigorous convergence analysis.

The resource allocation problem, a special case of the EMO problem, was studied for problems with smooth objective functions in \cite{arXiv:YHL}.
An intuitive extension of the continuous-time algorithm given in \cite{arXiv:YHL} to nonsmooth EMO cases may be written as
\begin{eqnarray}
\label{111}
\begin{cases}
\dot x_i(t) &\in \big{\{}p:p=P_{\Omega_i} [x_i(t)-g_i(x_i(t))+W_i^{\rm T}\lambda_i(t)]-x_i(t),\, g_i(x_i(t))\in \partial{f^i}(x_i(t))\big{\}},\\
\dot \lambda_i(t) &= d_i-W_i x_i(t)-\sum_{j=1}^{n}a_{i,j}(\lambda_i(t)-\lambda_j(t))-\sum_{j=1}^{n}a_{i,j}(z_i(t)-z_j(t)), \\
\dot z_i(t) &= \sum_{j=1}^{n}a_{i,j}(\lambda_i(t)-\lambda_j(t)),
\end{cases}
\end{eqnarray}
where $t\geq 0,\,i\in\{1,\ldots,n\},\,x_i(0)=x_{i0}\in\Omega_i\subset \mathbb R^{q_i}$, $\lambda_i(0)=\lambda_{i0}\in\mathbb R^m$, $z_i(0)=z_{i0}\in\mathbb R^m$, $\sum_{i=1}^n d_i=d_0$, and $a_{i,j}$ is the $(i,j)$th element of the   adjacency matrix of graph $\mathcal G$.
However, this algorithm involves the projection of subdifferential set (from $\partial{f^i}(x_i)$ to $\Omega_i$), which makes its convergence analysis very hard in the nonsmooth case. To overcome the technical challenges, we propose two different ideas to construct effective algorithms for the EMO problem in the following two subsections.

\subsection{Distributed Projected Output Feedback  Algorithm (DPOFA)}
\label{project_output}

The first idea is to use an auxiliary variable to avoid the projection of subdifferential set in the algorithm for EMO problem (\ref{optimization_problem}).
In other words, we propose a distributed  algorithm based on projected output feedbacks,  and the projected output feedback of the auxiliary variable is adopted to track the  optimal solution.
To be strict, we propose the continuous-time algorithm of agent $i$ as follows:
\begin{eqnarray}\label{distributed feedback_v2_x}
\begin{cases}
\dot y_i(t) &\in \big{\{}p:p=-y_i(t) + x_i(t)-g_i(x_i(t))+W_i^{\rm T}\lambda_i(t),\, g_i(x_i(t))\in \partial{f^i}(x_i(t))\big{\}},\\
\dot \lambda_i(t) &= d_i-W_i x_i(t)-\sum_{j=1}^{n}a_{i,j}(\lambda_i(t)-\lambda_j(t))-\sum_{j=1}^{n}a_{i,j}(z_i(t)-z_j(t)), \\
\dot z_i(t) &= \sum_{j=1}^{n}a_{i,j}(\lambda_i(t)-\lambda_j(t)),\\
x_i(t)&= P_{\Omega_i}(y_i(t)),
\end{cases}
\end{eqnarray}
with the auxiliary variable $y_i$(t) for $t\geq 0,\,y_i(0)=y_{i0}\in\mathbb R^{q_i},\,i\in\{1,\ldots,n\}$ and other notations are kept the same as those for (\ref{111}).
%where $\lambda_i(0)=\lambda_{i0}\in\mathbb R^m$, $z_i(0)=z_{i0}\in\mathbb R^m$, $\sum_{i=1}^n d_i=d_0$, and $a_{i,j}$ is the $(i,j)$th element of the   adjacency matrix of graph $\mathcal G$.
The term $x_i(t)= P_{\Omega_i}(y_i(t))$ is viewed as ``projected output feedback", which is inspired by \cite{LW:MAY2013}.  In this way, we avoid the technical difficulties resulting from the projection of the subdifferential.

Let $x \triangleq [x_1^{\rm T},\ldots,x_n^{\rm T}]^{\rm T}\in\Omega\subset\mathbb R^{\sum_{i=1}^{n}q_i}$, $y \triangleq [y_1^{\rm T},\ldots,y_n^{\rm T}]^{\rm T}\in\mathbb R^{\sum_{i=1}^{n}q_i}$, $ \lambda \triangleq [\lambda_1^{\rm T},\ldots,\lambda_n^{\rm T}]^{\rm T}\in\mathbb R^{nm}$, $d \triangleq [d_1^{\rm T},\ldots,d_n^{\rm T}]^{\rm T}\in\mathbb R^{nm}$, and $z \triangleq [z_1^{\rm T},\ldots,z_n^{\rm T}]^{\rm T}\in\mathbb R^{nm}$, where $\Omega\triangleq\prod_{i=1}^n \Omega_i$. Let $W = [W_1,\ldots,W_n]\in\mathbb{R}^{m\times \sum_{i=1}^{n}q_i}$ and $\overline{W} = \mathrm{diag}\{W_1,\ldots,W_n\}\in\mathbb{R}^{nm\times \sum_{i=1}^{n}q_i}$. Then (\ref{distributed feedback_v2_x}) can be written in a compact form
\begin{eqnarray}
\begin{bmatrix}\dot { y}(t) \\ \dot {\lambda}(t) \\ \dot z(t)\end{bmatrix}
 & \in & \mathcal F({y}(t),\lambda (t),z(t)),\label{feedback_comp_observer}\\
x(t) &=& P_{\Omega}(y(t)),\label{feedback_comp_proj}\\
\mathcal  F({y},\lambda,z )&\triangleq &\Bigg{\{}\begin{bmatrix} -y+x  -  g({x})+\overline{W}^{\rm T}\lambda\\
d-\overline{W}x-{L}{\lambda}-{L}{z}\\
{L}{\lambda}\end{bmatrix}: g({ x})\in\partial {f}( {x} ),\,x= P_{\Omega}(y)\Bigg{\}},\label{output}
\end{eqnarray}
where ${y}(0)={y}_0$, $\lambda (0)=\lambda _0$, ${z}(0)={z}_0$,
${L}=L_n\otimes I_m\in\mathbb R^{nm\times nm}$, and $L_n\in\mathbb R^{n\times n}$ is the  Laplacian matrix of graph $\mathcal G$.

\begin{remark}
In this algorithm, $x(t)=P_{\Omega}(y(t))$ is used to estimate the optimal solution of the EMO problem. Moreover, $x(t)$ stays in the constraint set $\Omega$ for every $t\geq 0$, though $y(t)$ may be out of $\Omega$.
Later, we will show that this algorithm can solve the EMO problem with nonsmooth objective functions and private constraints.
\end{remark}

\subsection{Distributed Derivative Feedback  Algorithm (DDFA)}\label{sec:alg_proj}
The second idea to facilitate the convergence analysis is to make a copy of the projection term by using the derivative feedback so as to  eliminate the ``trouble" caused by the projection term somehow.
To be strict, we propose the algorithm for the EMO problem (\ref{optimization_problem}) as follows:
\begin{eqnarray}\label{distributed project_v2_x}
\begin{cases}
\dot x_i(t) &\in \big{\{}p:p=P_{\Omega_i} [x_i(t)-g_i(x_i(t))+W_i^{\rm T}\lambda_i(t)]-x_i(t),\, g_i(x_i(t))\in \partial{f^i}(x_i(t))\big{\}},\\
\dot \lambda_i(t) &= d_i-W_i x_i(t)-\sum_{j=1}^{n}a_{i,j}(\lambda_i(t)-\lambda_j(t))-\sum_{j=1}^{n}a_{i,j}(z_i(t)-z_j(t))-W_i \dot x_i(t), \\
\dot z_i(t) &= \sum_{j=1}^{n}a_{i,j}(\lambda_i(t)-\lambda_j(t)),
\end{cases}
\end{eqnarray}
where all the  notations remain the same as  in (\ref{111}). %$t\geq 0,\,i\in\{1,\ldots,n\},\,x_i(0)=x_{i0}\in\Omega_i\subset \mathbb R^{q_i}$, $\lambda_i(0)=\lambda_{i0}\in\mathbb R^m$, $z_i(0)=z_{i0}\in\mathbb R^m$, $\sum_{i=1}^n d_i=d_0$, and $a_{i,j}$ is the $(i,j)$th element of the   adjacency matrix of graph $\mathcal G$.
Note that there is a derivative term $\dot x_i(t)$, viewed as ``derivative feedback", on the right-hand side of the second equation. The  derivative term can be found to be effective in the cancellation of the ``trouble" term $P_{\Omega_i} [x_i(t)-g_i(x_i(t))+W_i^{\rm T}\lambda_i(t)]$ in the analysis as demonstrated later.

Denote $x,\; \lambda, \: d,\; z,\; W$, and $\overline{W}$ as in \eqref{output}. Then (\ref{distributed project_v2_x}) can be written in a compact form
\begin{eqnarray}
\begin{bmatrix}
\dot { x}(t) \\ \dot {\lambda}(t) \\ \dot z(t)
\end{bmatrix}
 & \in & \mathcal F({x}(t),\lambda (t),z(t)),\label{project_comp_observer} \\
\mathcal  F({x},\lambda,z )
&\triangleq &\Bigg{\{}
\begin{bmatrix} p\\
d-\overline{W}x-{L}{\lambda}-{L}{z}-\overline{W}p\\
{L}{\lambda}\end{bmatrix}:\nonumber\\
&& p=P_{\Omega} [x-g(x)+\overline{W}^{\rm T}\lambda]-x,\,g({x})\in\partial {f}( {x} )\Bigg{\}},
\end{eqnarray}
where ${x}(0)={x}_0\in\Omega$, $\lambda (0)=\lambda _0\in\mathbb R^{nm}$, ${z}(0)={z}_0\in\mathbb R^{nm}$, ${L}=L_n\otimes I_m\in\mathbb R^{nm\times nm}$, and $L_n\in\mathbb R^{n\times n}$ is the  Laplacian matrix of graph $\mathcal G$.

\begin{assumption}\label{ASS_Exist}
For  algorithm (\ref{distributed project_v2_x}), or equivalently, algorithm (\ref{project_comp_observer}), the set-valued map  $\mathcal F({x},\lambda ,z)$ is with convex values for all $(x,\lambda,z)\in\Omega\times \mathbb{R}^{nm}\times \mathbb{R}^{nm}$.
\end{assumption}

Note that Assumption \ref{ASS_Exist} is given to guarantee the existence of solutions to algorithm \eqref{distributed project_v2_x}. In fact, in many situations, it can be satisfied. For example, it holds if $x_i\in\mathbb R$, or if both $\Omega_i$ and $\partial f^i_i(x_i)$ are ``boxes" for all $i\in\{1,\ldots,n\}$.
%\end{remark}

%We have the following observations for Algorithms \eqref{distributed feedback_v2_x} and \eqref{distributed project_v2_x}:
Besides the different techniques in algorithms \eqref{distributed feedback_v2_x} and \eqref{distributed project_v2_x}, the proposed algorithms are observed to  be different in the following aspects:
\begin{itemize}
  \item The application situations may be different: Although the convergence of both algorithms is based on the strict convexity assumption of the objective functions, algorithm \eqref{distributed project_v2_x} can also solve the EMO problem with only convex objective functions (which may have a continuum set of optimal solutions) when the objective functions are differentiable (see {Corollary} \ref{theorem_smooth_convergence} in Section \ref{main_result}).
  %\item Algorithm \eqref{distributed feedback_v2_x} changes variable $y$ and uses its projection $x=P_{\Omega}(y)$ to indirectly track the optimal solution. While algorithm \eqref{distributed project_v2_x} directly changes $x(t)\in\Omega$ to estimate the optimal solution, which shows faster  response speed of $x(t)$ than algorithm \eqref{distributed feedback_v2_x} under some conditions (see simulation results in Section \ref{sec_NS}).
  \item The dynamic performances may be different: Because algorithm \eqref{distributed project_v2_x} directly changes $x(t)\in\Omega$ to estimate the optimal solution, it may show faster  response speed of $x(t)$ than that of algorithm \eqref{distributed feedback_v2_x} (see simulation results in Section \ref{sec_NS}).
\end{itemize}

Furthermore, algorithms \eqref{distributed feedback_v2_x} and \eqref{distributed project_v2_x} are essentially different from existing ones. Compared with the algorithm in \cite{CC:network:2015}, our algorithms need not exchange information of subgradients among the agents.  Unlike the algorithm given in \cite{YHL:SCL:2015}, ours use different techniques (i.e., the projected output feedback in \eqref{distributed feedback_v2_x} or derivative feedback in \eqref{distributed project_v2_x}) to estimate the optimal solution. %Compared with \cite{NC:CDC:2015}, our algorithms avoid the orthogonal set projection. %, which lowers the computational cost.
Moreover, our algorithms have two advantages compared with previous methods.
\begin{itemize}
  \item Agent $i$ of the proposed algorithms knows $W_i$, which is composed of a subset of columns in $W$. This is different from existing results with assuming that each agent knows a subset of rows of the equality constraints \cite{MLM:TAC:2015}. %,LW:TAC:2015,GC:2014}.
      If $n$ is a sufficiently large number and $m $ is relatively small, the proposed designs make the computation load at  each node relatively small compared with previous algorithms in \cite{MLM:TAC:2015}. %,LW:TAC:2015,GC:2014}.
  \item  Agent $i$ in the proposed algorithms exchanges information of $\lambda _i\in\mathbb R^{m}$ and $z_i\in\mathbb R^{m}$ with its neighbors. Compared with algorithms which require exchanging $x_i\in\mathbb R^{q_i}$, this design can greatly reduce the communication cost when $q_i$ is much larger than $m$ for $i\in\{1,\ldots,n\}$ and the information of $x_i\in\mathbb R^{q_i}$ is kept  confidential.
\end{itemize}

\section{Convergence Analysis}
\label{main_result}
In this section, we use the stability analysis of differential inclusions to prove the correctness and convergence of our proposed algorithms.

\subsection{Convergence Analysis of DPOFA}%Projected Output Feedback Algorithm}

Consider algorithm \eqref{distributed feedback_v2_x} (or \eqref{feedback_comp_observer}). Let $(y^{*},\lambda^{*},z^{*})\in\mathbb{R}^{\sum_{i=1}^{n}q_i}\times \mathbb{R}^{nm}\times \mathbb{R}^{nm}$ be an equilibrium of \eqref{distributed feedback_v2_x}. Then
\begin{subequations}\label{distributed_feedback_v2}
\begin{align}
% \nonumber % Remove numbering (before each equation)
  0_{\sum_{i=1}^{n}q_i} &\in \{p: p=-y^* + x^*  -  g({x^*})+\overline{W}^{\rm T} \lambda^* ,\,x^*=P_{\Omega}(y^*),\,g({ x^*})\in\partial {f}( {x^*} )\}, \label{distributed feedback_v2_x1}\\
  0_{nm} &= d-\overline{W}x^*-Lz^* ,\quad x^*=P_{\Omega}(y^*),\label{distributed feedback_v2_x2}\\
  0_{nm} &= L\lambda^*.\label{distributed feedback_v2_x3}
\end{align}
\end{subequations}

The following result reveals the relationship between the equilibrium points of algorithm \eqref{distributed feedback_v2_x} and the solutions of problem \eqref{optimization_problem}.

\begin{theorem}\label{equivalent_conditions}
If $(y^{*},\lambda^{*},z^{*})\in\mathbb{R}^{\sum_{i=1}^{n}q_i}\times \mathbb{R}^{nm}\times \mathbb{R}^{nm}$ is an equilibrium of \eqref{distributed feedback_v2_x}, then $x^{*}=P_{\Omega}(y^{*})$ is a solution to problem \eqref{optimization_problem}. Conversely, if $x^{*}\in\Omega$ is a solution to problem \eqref{optimization_problem}, then there exists $(y^{*}, \lambda^{*},z^{*})\in\mathbb{R}^{\sum_{i=1}^{n}q_i}\times \mathbb{R}^{nm}\times \mathbb{R}^{nm}$ such that $(y^{*},\lambda^{*},z^{*})$ is an equilibrium of \eqref{distributed feedback_v2_x} with $x^{*}=P_{\Omega}(y^{*})$.
\end{theorem}

\begin{IEEEproof}
($i$) Suppose $(y^{*},\lambda^{*},z^{*})\in\mathbb{R}^{\sum_{i=1}^{n}q_i}\times \mathbb{R}^{nm}\times \mathbb{R}^{nm}$ is an equilibrium of \eqref{distributed feedback_v2_x}. Left-multiply both sides of \eqref{distributed feedback_v2_x2} by $1_n^{\rm T}\otimes I_m$, it follows that
\begin{eqnarray}
1_n^{\rm T}\otimes I_m(d-\overline{W}x^*-Lz^*)&=&\sum_{i=1}^{n}(d_i-W_i x_i^*)-(1_n^{\rm T}\otimes I_m )Lz^*\nonumber\\
&=&d_0-Wx^*-(1_n^{\rm T}\otimes I_m) Lz^*=0_m.\label{d_x1}
\end{eqnarray}
It follows from the properties of Kronecker product and $L_n^{\rm T}{1}_n=0_n$ that
\begin{eqnarray}\label{d_x2}
(1_n^{\rm T}\otimes I_m )L=(1_n^{\rm T}\otimes I_m )(L_n\otimes I_m)=(1_n^{\rm T}L_n)\otimes(I_m)=0_{m\times nm}.
\end{eqnarray}
In the light of  (\ref{d_x1}) and (\ref{d_x2}), \eqref{opt_equi_2} holds.

Next, it follows from \eqref{distributed feedback_v2_x3} that there exists $\lambda_0^*\in\mathbb{R}^m$ such that $\lambda^*=1_n\otimes\lambda_0^*$.
By taking into consideration  \eqref{distributed feedback_v2_x1} and $\lambda^*=1_n\otimes\lambda_0^*$, there exists $g(x^*)\in\partial f(x^*)$ such that $x^*  -  g({x^*})+\overline{W}^{\rm T}(1_n\otimes\lambda_0^*)=x^*  -  g({x^*})+{W}^{\rm T}\lambda_0^*=y^*$. Since $x^*=P_{\Omega}(y^*)$, it follows that \eqref{opt_equi_1} holds.

By virtue of  \eqref{opt_equi_1}, \eqref{opt_equi_2} and Lemma \ref{equiv_cond}, $x^*$ is the solution to problem \eqref{optimization_problem}.

($ii$) Conversely, suppose $x^*$ is the solution to problem \eqref{optimization_problem}. According to Lemma \ref{equiv_cond}, there exist $\lambda^*_0\in\mathbb R^{q}$ and $g(x^*)\in\partial f(x^*)$ such that \eqref{opt_equi_1} and \eqref{opt_equi_2} hold. Define $\lambda^*=1_n\otimes\lambda_0^*$. As a result, \eqref{distributed feedback_v2_x3} holds.
%It follows that $W^{\rm T}\lambda_0^*=\overline{W}^{\rm T}(1_n\otimes\lambda_0^*)=\overline{W}^{\rm T}\lambda^*$ and, hence, \eqref{distributed feedback_v2_x3} holds.

Take any $v\in\mathbb{R}^m$ and let $\mathbf v = 1_n\otimes v$. Since \eqref{opt_equi_2} holds, $(d-\overline{W}x^*)^{\rm T}\mathbf v=(\sum_{i=1}^{n}(d_i-W_ix_i^*))^{\rm T} v=(d_0-Wx^*)^{\rm T}v=0$. Due to the properties of Kronecker product and $L_n {1}_n=0_n$, $L\mathbf v=(L_n\otimes I_m)(1_n\otimes
v)=(L_n1_n)\otimes(I_mv)=0_n\otimes v=0_{nm}$ and, hence, $\mathbf v\in\ker(L)$.
Note that $\ker (L)$ and $\mathrm{range} (L)$ form an orthogonal decomposition of $\mathbb R^{nm}$ by the fundamental theorem of linear algebra \cite{Strang:1993}. It follows from $(d-\overline{W} x^*)^{\rm T}\mathbf v=0$ and $\mathbf v\in\ker(L)$ that $d-\overline{W}x^*\in\mathrm{range} (L)$.
Hence, there exists $z^*\in\mathbb R^{nm}$ such that \eqref{distributed feedback_v2_x2} holds.

Because $W^{\rm T}\lambda_0^*=\overline{W}^{\rm T}(1_n\otimes\lambda_0^*)=\overline{W}^{\rm T}\lambda^*$, \eqref{opt_equi_1} implies $x^*
%=P_{\Omega}(x^*  -  g({x^*})+{W}^{\rm T} \lambda_0^*)
=P_{\Omega}(x^*  -  g({x^*})+\overline{W}^{\rm T} (1_n\otimes\lambda_0^*))=P_{\Omega}(x^*  -  g({x^*})+\overline{W}^{\rm T}\lambda^*)$ for some $g({x^*})\in\partial f(x^*)$. Choose $y^* = x^*  -  g({x^*})+\overline{W}^{\rm T}\lambda^*$. \eqref{distributed feedback_v2_x1} holds.

To sum up, if $x^{*}\in\Omega$ is a solution to problem \eqref{optimization_problem}, there exists $(y^{*}, \lambda^{*},z^{*})\in\mathbb{R}^{\sum_{i=1}^{n}q_i}\times \mathbb{R}^{nm}\times \mathbb{R}^{nm}$ such that \eqref{distributed_feedback_v2} holds and $x^{*}=P_{\Omega}(y^{*})$. Hence, $(y^{*},\lambda^{*},z^{*})$ is an equilibrium of \eqref{distributed feedback_v2_x} with $x^{*}=P_{\Omega}(y^{*})$.
\end{IEEEproof}

Let $x^{*}$ be the solution to problem \eqref{optimization_problem}. It follows from Theorem \ref{equivalent_conditions} that there exist $y^{*}$, $\lambda^{*}$ and $z^{*}$ such that $(y^{*},\lambda^{*},z^{*})$ is an equilibrium of \eqref{distributed feedback_v2_x} with $x^{*}=P_{\Omega}(y^{*})$. Define the function
\begin{eqnarray}\label{V}
V(y,\lambda,z)\triangleq\frac{1}{2}(\|y-P_{\Omega}(y^*)\|^2-\|y-P_{\Omega}(y)\|^2)+\frac{1}{2}\|\lambda-\lambda^*\|^2+\frac{1}{2}\|z-z^*\|^2. \end{eqnarray}

\begin{lemma}\label{differential_eqn}
Consider  algorithm (\ref{distributed feedback_v2_x}). Under Assumption \ref{Assumption} with $V(y,\lambda,z)$ defined in \eqref{V},
if $a\in\mathcal L_{{F}} V(y,\lambda,z)$,  then there exist $g({x})\in\partial {f}({x})$ and $g({x^*})\in\partial {f}({x^*})$ with $x=P_{\Omega}(y)$ and $x^*=P_{\Omega}(y^*)$ such that $a\leq  -(x-x^*)^{\rm T}(g({x})-g({x^*}))-\lambda^{\rm T}L\lambda\leq 0$.
\end{lemma}
\begin{IEEEproof}
It follows from Lemma \ref{projection_diff} that the gradient of $V(y,\lambda,z)$ with respect to $y$ is $\nabla_y V(y,\lambda,z)=x-x^*$ where $x=P_{\Omega}(y)$ and $x^*=P_{\Omega}(y^*)$. The gradients of $V(y,\lambda,z)$ with respect to $\lambda$ and $z$ are $\nabla_{\lambda} V(y,\lambda,z)=\lambda-\lambda^*$ and $\nabla_{z} V(y,\lambda,z)=z-z^*$, respectively.

The function $V(y,\lambda,z)$ along the trajectories of (\ref{distributed feedback_v2_x})  satisfies that
\begin{eqnarray*}
% \nonumber % Remove numbering (before each equation)
  \mathcal L_{\mathcal F}V(y,\lambda,z) &=& \Big{\{}a\in\mathbb{R}:a=\nabla_y V(y,\lambda,z)^{\rm T}(-y+x  -  g({x})+\overline{W}^{\rm T}\lambda)\\
  &&+\nabla_{\lambda} V(y,\lambda,z)^{\rm T}(d-\overline{W} x-{L}{\lambda}-{L}{z})+\nabla_{z} V(y,\lambda,z)^{\rm T}{L}{\lambda},\\
  &&\,g({ x})\in\partial {f}( {x} ),\,x= P_{\Omega}(y)\Big{\}}.
\end{eqnarray*}

Suppose $a\in\mathcal L_{\mathcal F}V(y,\lambda,z)$,  then there is $g({ x})\in\partial {f}( {x} )$ such that
\begin{eqnarray}\label{diff_V}
% \nonumber % Remove numbering (before each equation)
  a = (x-x^*)^{\rm T}(-y+x  -  g({x})+\overline{W}^{\rm T} \lambda)+(\lambda-\lambda^*)^{\rm T}(d-\overline{W} x-{L}{\lambda}-{L}{z})+(z-z^*)^{\rm T}{L}{\lambda},
\end{eqnarray}
where $x=P_{\Omega}(y)$.

Because $(y^{*},\lambda^{*},z^{*})$ is an equilibrium of \eqref{distributed feedback_v2_x} with $x^{*}=P_{\Omega}(y^{*})$, there exists $g(x^*)\in\partial f(x^*)$ such that
\begin{eqnarray}\label{equi_c}
\begin{cases}
0_{nm} &= L\lambda^{*}\\
d &= \overline{W} x^{*}+{L}{z^*},\\
0_{\sum_{i=1}^n q_i} &= -y^*+x^*  -  g({x^*})+\overline{W}^{\rm T} \lambda^*.
\end{cases}
\end{eqnarray}

It follows from \eqref{diff_V} and \eqref{equi_c} that
\begin{eqnarray}\label{diff_V2}
% \nonumber % Remove numbering (before each equation)
  a &=& (x-x^*)^{\rm T}\big{[}(-y+x  -  g({x})+\overline{W}^{\rm T}\lambda)- (-y^*+x^*  -  g({x^*})+\overline{W}^{\rm T}\lambda^*)\big{]}\nonumber\\
  &&+(\lambda-\lambda^*)^{\rm T}(\overline{W} x^{*}+{L}{z^*}-\overline{W} x-{L}{\lambda}-{L}{z}) +(z-z^*)^{\rm T}{L}{\lambda}\nonumber\\
  &=& -(x-x^*)^{\rm T}(y-y^*)+\|x-x^*\|^2-(x-x^*)^{\rm T}(g({x})-g({x^*}))+(x-x^*)^{\rm T}\overline{W}^{\rm T}(\lambda-\lambda^*)\nonumber\\
  &&-(x-x^*)^{\rm T}\overline{W}^{\rm T}(\lambda-\lambda^*)-\lambda^{\rm T}L\lambda-\lambda^{\rm T}L(z-z^*)+(z-z^*)^{\rm T}{L}{\lambda}\nonumber\\
  &=& -(x-x^*)^{\rm T}(y-y^*)+\|x-x^*\|^2-(x-x^*)^{\rm T}(g({x})-g({x^*}))-\lambda^{\rm T}L\lambda.
\end{eqnarray}

Since $x=P_{\Omega}(y)$ and $x^{*}=P_{\Omega}(y^{*})$, we obtain from Lemma \ref{lemma_ineq} that $-(x-x^*)^{\rm T}(y-y^*)+\|x-x^*\|^2\leq 0$. Hence,
\begin{eqnarray*}
% \nonumber % Remove numbering (before each equation)
  a \leq -(x-x^*)^{\rm T}(g({x})-g({x^*}))-\lambda^{\rm T}L\lambda.
\end{eqnarray*}

The convexity of $f$ implies that $(x-x^*)^{\rm T}(g({x})-g({x^*}))\geq 0$. In addition, $L=L_n\otimes I_q\geq 0$ since $L_n\geq 0$. Hence, $a\leq -(x-x^*)^{\rm T}(g({x})-g({x^*}))-\lambda^{\rm T}L\lambda\leq 0$.
\end{IEEEproof}

The following result shows the correctness of the proposed algorithm.
\begin{theorem}\label{theorem_convergence}
Consider  algorithm (\ref{distributed feedback_v2_x}). If Assumption \ref{Assumption} holds, then
\begin{itemize}
\item[($i$)] every solution $({y}(t),x(t),\lambda(t),z(t))$ is bounded;
\item[($ii$)] for every solution, ${x}(t)$ converges to the optimal solution to problem \eqref{optimization_problem}.
\end{itemize}
\end{theorem}

\begin{IEEEproof}
($i$) Let $V(y,\lambda,z)$ be as defined in \eqref{V}. It follows from Lemma \ref{differential_eqn} that %
\begin{eqnarray}\label{LV_negg}
  \max\mathcal L_{\mathcal F}V(y,\lambda,z)\leq \max \{-(x-x^*)^{\rm T}(g({x})-g({x^*}))-\lambda^{\rm T}L\lambda:g(x)\in\partial f(x)\}\leq 0.
\end{eqnarray}

Note that $V(y,\lambda,z)\geq \frac{1}{2}\|x-x^*\|^2+\frac{1}{2}\|\lambda-\lambda^*\|^2+\frac{1}{2}\|z-z^*\|^2$ in view of Lemma \ref{projection_diff}. It follows from \eqref{LV_negg} that $(x(t),\lambda(t),z(t)),\,t\geq 0$ is bounded. Because $\partial f(x)$ is compact, there exists $m=m(y_0,\lambda_0,z_0)>0$ such that
\begin{eqnarray}\label{m}
\|x(t)  -  g({x(t)})+\overline{W}^{\rm T}\lambda(t)\|<m,
\end{eqnarray}
for all $g({x(t)})\in\partial f(x(t))$ and all $t\geq 0.$

Define $X:\mathbb{R}^{\sum_{i=1}^{n}q_i}\rightarrow\mathbb R$ by $X(y)=\frac{1}{2}\|y\|^2$. The function $X(y)$ along the trajectories of (\ref{distributed feedback_v2_x})  satisfies that
\begin{eqnarray*}
% \nonumber % Remove numbering (before each equation)
  \mathcal{L}_{\mathcal F} X(y) =\{y^{\rm T}(-y+x  -  g({x})+\overline{W}^{\rm T}\lambda):g(x)\in\partial f(x)\}.
\end{eqnarray*}
Note that $y^{\rm T}(t)(-y(t)+x(t)  -  g({x(t)})+\lambda(t))\leq -\|y(t)\|^2+m\|y(t)\|$, where $t\geq 0$, $m$ is defined by \eqref{m} and $g(x(t))\in\partial f(x(t))$. Hence,
$$\max\mathcal{L}_{\mathcal F} X(y(t))\leq -\|y(t)\|^2+m\|y(t)\|=-2X(y(t))+m\sqrt{2X(y(t))}.$$
It can be easily verified that $X(y(t)),\,t\geq 0,$ is bounded, so is $y(t),\,t\geq 0$.

Part ($i$) is thus proved.

($ii$) Let $\mathcal{R}=\big{\{} (y,\lambda,z)\in\mathbb{R}^{\sum_{i=1}^{n}q_i}\times \mathbb{R}^{nm}\times \mathbb{R}^{nm}: 0\in\mathcal{L}_{\mathcal F}V(y,\lambda,z) \big{\}}\subset \big{\{} (y,\lambda,z)\in\mathbb{R}^{\sum_{i=1}^{n}q_i}\times \mathbb{R}^{nm}\times \mathbb{R}^{nm}: 0=\min_{g({x})\in\partial f(x),\,g({x^*})\in\partial f(x^*)}(x-x^*)^{\rm T}(g({x})-g({x^*})),\,L\lambda=0,\,x=P_{\Omega}(y),\, x^{*}=P_{\Omega}(y^{*})\big{\}}$. Note that $(x-x^*)^{\rm T}(g({x})-g({x^*}))>0$ if $x\not=x^*$ because of the strict convexity assumption of $f$ and Lemma \ref{convex_strict}. Hence, $\mathcal{R}\subset\big{\{} (y,\lambda,z)\in\mathbb{R}^{\sum_{i=1}^nq_i}\times \mathbb{R}^{nm}\times \mathbb{R}^{nm}: x=P_{\Omega}(y)=x^* ,\,L\lambda=0\big{\}}$.

Let $\mathcal M$ be the largest weakly invariant subset of $\overline{\mathcal R}$. It follows from Lemma \ref{nonsmooth_invariance} that $(y(t),\lambda(t),z(t))\rightarrow\mathcal M$ as $t\rightarrow \infty$. Hence, $x(t)\rightarrow x^*$ as $t\rightarrow\infty$.

Part ($ii$) is thus proved.
\end{IEEEproof}

%\begin{remark}
%The proposed algorithm is essentially different from existing results given in \cite{XJB:2004,YHL:SCL:2015,arXiv:YHL,CC:network:2015}, which leads to a different technique in the convergence proof. Our result is inspired by the result in \cite{LW:MAY2013}. Because our algorithm is distributed and considers nonsmooth objective functions and more general equality constraints, our proposed algorithm is new.
%\end{remark}

\subsection{Convergence Analysis of DDFA}

Consider algorithm \eqref{distributed project_v2_x} (or \eqref{project_comp_observer}). $(x^{*},\lambda^{*},z^{*})\in\Omega\times \mathbb{R}^{nm}\times \mathbb{R}^{nm}$ is an equilibrium of \eqref{distributed project_v2_x}  if and only if there exists $g({ x^*})\in\partial {f}( {x^*} )$ such that
\begin{subequations}\label{distributed_project_v2}
\begin{align}
% \nonumber % Remove numbering (before each equation)
  0_{\sum_{i=1}^{n}q_i} &=P_{\Omega} [x^*  -  g({x^*})+\overline{W}^{\rm T} \lambda^*]- x^*, \label{distributed project_v2_x1}\\
  0_{nm} &= d-\overline{W}x^*-Lz^*,\label{distributed project_v2_x2}\\
  0_{nm} &= L\lambda^*.\label{distributed project_v2_x3}
\end{align}
\end{subequations}

\begin{theorem}\label{equivalent_project_conditions}
If $(x^{*},\lambda^{*},z^{*})\in\Omega\times \mathbb{R}^{nm}\times \mathbb{R}^{nm}$ is an equilibrium of \eqref{distributed project_v2_x}, then $x^{*}$ is a solution to problem \eqref{optimization_problem}. Conversely, if $x^{*}\in\Omega$ is a solution to problem \eqref{optimization_problem}, then there exists $\lambda^{*}\in\mathbb{R}^{nm}$ and $z^{*}\in \mathbb{R}^{nm}$ such that $(x^{*},\lambda^{*},z^{*})$ is an equilibrium of \eqref{distributed project_v2_x}.
\end{theorem}

The proof is similar to that of Theorem \ref{equivalent_conditions} and, hence, is omitted.

Suppose $(x^{*},\lambda^{*},z^{*})\in\Omega\times \mathbb{R}^{nm}\times \mathbb{R}^{nm}$ is an equilibrium of \eqref{distributed project_v2_x}. Define the  function
\begin{eqnarray}\label{Lya_func}
V(x,\lambda,z)=f(x)-f(x^*)+(\lambda^*)^{\rm T}(d-\overline{W}x) +\frac{1}{2}\|x-x^*\|^2+\frac{1}{2}\|\lambda-\lambda^*\|^2+\frac{1}{2}\|z-z^*\|^2.
\end{eqnarray}

\begin{lemma}\label{V_positive}
Let function $V(x,\lambda,z)$ be as defined in \eqref{Lya_func} and Assumption \ref{Assumption} hold. For all $(x,\lambda,z)\in \Omega\times \mathbb{R}^{nm}\times \mathbb{R}^{nm}$, $V(x,\lambda,z)$ is positive definite, $V(x,\lambda,z)=0$ if and only if $(x,\lambda,z)=(x^{*},\lambda^{*},z^{*})$, and $V(x,\lambda,z)\rightarrow\infty$ as $(x,\lambda,z)\rightarrow\infty$.
\end{lemma}

\begin{IEEEproof}
By \eqref{distributed project_v2_x3}, there is $\lambda_0^*\in\mathbb{R}^m$ such that $\lambda^*=1_n\otimes\lambda_0^*$.
It can be easily verified that $(\lambda^*)^{\rm T}(d-\overline{W}x)=(\lambda_0^*)^{\rm T}(d_0-{W}x)$.
It follows from \eqref{opt_equi_2}   that
\begin{eqnarray}\label{equi_2_eq}
f(x)-f(x^*)+(\lambda^*)^{\rm T}(d-\overline{W}x)=f(x)-f(x^*)+(\lambda^*_0)^{\rm T}{W}(x^*-x).
\end{eqnarray}
By \eqref{Lya_func} and \eqref{equi_2_eq}, it is straightforward that $V(x^*,\lambda^*,z^*)=0$.

Because $f(x)$ is convex, $f(x)-f(x^*)\geq g(x^*)^{\rm T}(x-x^*) $ for all $x\in\Omega$ and $g(x^*)\in\partial f(x^*)$. According to \eqref{equi_1},
\begin{eqnarray}\label{equi_1_eq}
(g(x^*)-W^{\rm T}\lambda^*_0)^{\rm T}(x-x^*)\geq 0
\end{eqnarray}
for all $x\in\Omega$, where $g(x^*)\in\partial f(x^*)$ is as chosen in \eqref{equi_1}.
It follows from \eqref{equi_2_eq} and \eqref{equi_1_eq} that
\begin{eqnarray}
f(x)-f(x^*)+(\lambda^*)^{\rm T}(d-\overline{W}x)\geq (g(x^*)-W^{\rm T}\lambda^*_0)^{\rm T}(x-x^*)\geq 0
\end{eqnarray}
for all $x\in\Omega$, where $g(x^*)\in\partial f(x^*)$ is as chosen in \eqref{equi_1}.

Hence,  for all $(x,\lambda,z)\in\Omega\times \mathbb{R}^{nm}\times \mathbb{R}^{nm}$, $V(x,\lambda,z)\geq \frac{1}{2}\|x-x^*\|^2+\frac{1}{2}\|\lambda-\lambda^*\|^2+\frac{1}{2}\|z-z^*\|^2$. Therefore, $V(x,\lambda,z)$ is positive definite, $V(x,\lambda,z)=0$ if and only if $(x,\lambda,z)=(x^{*},\lambda^{*},z^{*})$, and $V(x,\lambda,z)\rightarrow\infty$ as $(x,\lambda,z)\rightarrow\infty$ for all $(x,\lambda,z)\in\Omega\times \mathbb{R}^{nm}\times \mathbb{R}^{nm}$.
\end{IEEEproof}

\begin{lemma}\label{differential_project_eqn}
Consider algorithm (\ref{distributed project_v2_x}). Under Assumptions \ref{Assumption} and \ref{ASS_Exist} with $V(x,\lambda,z)$ defined in \eqref{Lya_func},
if $a\in\mathcal L_{{F}} V(x,\lambda,z)$,  then there exist $g({x})\in\partial {f}({x})$ and $g({x^*})\in\partial {f}({x^*})$ such that $a\leq -\|p\|^2 -(x-x^*)^{\rm T}(g({x})-g({x^*}))-\lambda^{\rm T}L\lambda\leq 0$, where $p=P_{\Omega} [x  -  g({x})+\overline{W}^{\rm T} \lambda]- x$.
\end{lemma}

\begin{IEEEproof}
The function $V(x,\lambda,z)$ along the trajectories of (\ref{distributed project_v2_x})  satisfies that
\begin{eqnarray*}
% \nonumber % Remove numbering (before each equation)
  \mathcal L_{\mathcal F}V(x,\lambda,z) &=& \Big{\{}a\in\mathbb{R}:a=(g(x)-\overline{W}^{\rm T}\lambda^{*}+x-x^*)^{\rm T}p\\
  &&+\nabla_{\lambda} V(x,\lambda,z)^{\rm T}(d-\overline{W} x-{L}{\lambda}-{L}{z}-\overline{W}p)+\nabla_{z} V(x,\lambda,z)^{\rm T}{L}{\lambda},\\
  &&\,g({ x})\in\partial {f}( {x} ),\,p=P_{\Omega} [x  -  g({x})+\overline{W}^{\rm T} \lambda]- x\Big{\}}.
\end{eqnarray*}

Suppose $a\in\mathcal L_{\mathcal F}V(x,\lambda,z)$,  then there exists $g({ x})\in\partial {f}( {x} )$ such that
\begin{eqnarray}\label{diff_V_project}
% \nonumber % Remove numbering (before each equation)
  a = (g(x)-\overline{W}^{\rm T}\lambda^{*}+x-x^*)^{\rm T}p+(\lambda-\lambda^*)^{\rm T}(d-\overline{W} x-{L}{\lambda}-{L}{z}-\overline{W}p)+(z-z^*)^{\rm T}{L}{\lambda},
\end{eqnarray}
where
\begin{eqnarray}\label{p_equ_P_x}
p=P_{\Omega} [x  -  g({x})+\overline{W}^{\rm T} \lambda]- x.
\end{eqnarray}

By \eqref{projection_inequality}, we represent \eqref{p_equ_P_x} in the form of a variational inequality
\begin{eqnarray*}
% \nonumber % Remove numbering (before each equation)
  \langle p+x-(x  -  g({x})+\overline{W}^{\rm T} \lambda) , p+x -\tilde{x}\rangle \leq 0,\quad \forall \tilde{x}\in\Omega.
\end{eqnarray*}
Choose $\tilde{x}=x^*$. Then,
\begin{eqnarray}\label{inq_p}
% \nonumber % Remove numbering (before each equation)
  (g(x)-\overline{W}^{\rm T}\lambda+x-x^*)^{\rm T}p\leq -\|p\|^2-(g(x)-\overline{W}^{\rm T}\lambda)^{\rm T}(x-x^*).
\end{eqnarray}

Since $(x^{*},\lambda^{*},z^{*})$ is an equilibrium of \eqref{distributed project_v2_x}, there is $g(x^*)\in\partial f(x^*)$ such that
\begin{eqnarray}\label{equi_p}
\begin{cases}
0_{nm} & = L\lambda^{*}\\
d & = \overline{W} x^{*}+{L}{z^*},\\
x^{*} & =P_{\Omega} [x^{*}  -  g({x}^{*})+\overline{W}^{\rm T} \lambda^{*}].
\end{cases}
\end{eqnarray}

It follows from \eqref{diff_V_project}, \eqref{inq_p} and \eqref{equi_p} that
\begin{eqnarray}\label{diff_V2_project}
% \nonumber % Remove numbering (before each equation)
  a &=& (g(x)-\overline{W}^{\rm T}\lambda+x-x^*)^{\rm T}p+(\lambda-\lambda^{*})^{\rm T}\overline{W}p\nonumber\\
  &&+(\lambda-\lambda^*)^{\rm T}(\overline{W} x^{*}+{L}{z^*}-\overline{W}^{\rm T} x-{L}{\lambda}-{L}{z}-\overline{W}p) +(z-z^*)^{\rm T}{L}{\lambda}\nonumber\\
  &\leq & -\|p\|^2-(g(x)-\overline{W}^{\rm T}\lambda)^{\rm T}(x-x^*)+(\lambda-\lambda^{*})^{\rm T}\overline{W}p-(x-x^*)^{\rm T}\overline{W}^{\rm T}(\lambda-\lambda^*)-\lambda^{\rm T}L\lambda
  \nonumber\\
  &&-\lambda^{\rm T}L(z-z^*)-(\lambda-\lambda^{*})^{\rm T}\overline{W}p+(z-z^*)^{\rm T}{L}{\lambda}\nonumber\\
  &=& -\|p\|^2-(g(x)-\overline{W}^{\rm T}\lambda)^{\rm T}(x-x^*)-(x-x^*)^{\rm T}\overline{W}^{\rm T}(\lambda-\lambda^*)-\lambda^{\rm T}L\lambda\nonumber\\
  &=& -\|p\|^2-(g(x)-g(x^*))^{\rm T}(x-x^*)-(g(x^*)-\overline{W}^{\rm T}\lambda^*)^{\rm T}(x-x^*)-\lambda^{\rm T}L\lambda.
\end{eqnarray}

Because $x^{*}  =P_{\Omega} [x^{*}  -  g({x}^{*})+\overline{W}^{\rm T} \lambda^{*}]$, $(g(x^*)-\overline{W}^{\rm T}\lambda^*)^{\rm T}(x-x^*)\geq 0$ for all $x\in\Omega$ followed by \eqref{projection_inequality}.
The convexity of $f$ implies that $(x-x^*)^{\rm T}(g({x})-g({x^*}))\geq 0$. In addition, $L=L_n\otimes I_q\geq 0$ since $L_n\geq 0$. Hence, $a\leq -\|p\|^2-(x-x^*)^{\rm T}(g({x})-g({x^*}))-\lambda^{\rm T}L\lambda\leq 0$.
\end{IEEEproof}

\begin{theorem}\label{theorem_project_convergence}
Consider  algorithm (\ref{distributed project_v2_x}). If Assumptions \ref{Assumption} and \ref{ASS_Exist} hold, then
\begin{itemize}
\item[($i$)] every solution $(x(t),\lambda(t),z(t))$ is bounded;
\item[($ii$)] for every solution, ${x}(t)$ converges to the optimal solution to problem \eqref{optimization_problem}.
\end{itemize}
\end{theorem}

\begin{IEEEproof}
($i$) Suppose $(x^{*},\lambda^{*},z^{*})\in\Omega\times \mathbb{R}^{nm}\times \mathbb{R}^{nm}$ is an equilibrium of \eqref{distributed project_v2_x}. Let function $V(x,\lambda,z)$ be as defined in \eqref{Lya_func}. It follows from Lemma \ref{differential_project_eqn} that
\begin{eqnarray*}
 \max \mathcal L_{\mathcal F}V (x,\lambda,z) &\leq &\sup \big{\{}a:a=-\|p\|^2 -(x-x^*)^{\rm T}(g({x})-g({x^*}))-\lambda^{\rm T}L\lambda,\\
 && \,g({ x})\in\partial {f}( {x} ),\,p=P_{\Omega} [x  -  g({x})+\overline{W}^{\rm T} \lambda]- x\big{\}}\leq 0.
\end{eqnarray*}
It follows from Lemma \ref{V_positive} that  $V(x,\lambda,z)$ is positive definite for all $(x,\lambda,z)\in\Omega\times \mathbb{R}^{nm}\times \mathbb{R}^{nm}$, $V(x,\lambda,z)=0$ if and only if $(x,\lambda,z)=(x^{*},\lambda^{*},z^{*})$, and $V(x,\lambda,z)\rightarrow\infty$ as $(x,\lambda,z)\rightarrow\infty$.  Hence, $(x(t),\lambda(t),z(t))$ is bounded for all $t\geq 0$.

($ii$) Let $\mathcal{R}=\big{\{} (x,\lambda,z)\in\Omega\times \mathbb{R}^{nm}\times \mathbb{R}^{nm}: 0\in\mathcal{L}_{\mathcal F}V(y,\lambda,z) \big{\}}\subset \big{\{} (x,\lambda,z)\in\Omega\times \mathbb{R}^{nm}\times \mathbb{R}^{nm}: \exists g({ x})\in\partial {f}( {x} ),\, L\lambda=0_{nq},(x-x^*)^{\rm T}(g({x})-g({x^*}))=0,\,\,0_{\sum_{i=1}^{n}q_i}=P_{\Omega} [x  -  g({x})+\overline{W}^{\rm T} \lambda]- x\big{\}}$. Let $\mathcal M$ be the largest weakly invariant subset of $\overline{\mathcal R}$.
It follows from Lemma \ref{nonsmooth_invariance} that $(x(t),\lambda(t),z(t))\rightarrow\mathcal M$ as $t\rightarrow \infty$. Note that $(x-x^*)^{\rm T}(g({x})-g({x^*}))>0$ if $x\not=x^*$ because of the strict convexity assumption of $f$ and Lemma \ref{convex_strict}. Hence, $x(t)\rightarrow x^*$ as $t\rightarrow\infty$.
\end{IEEEproof}

The following gives a convergence result when the objective functions of problem \eqref{optimization_problem} are differentiable. If the objective functions of problem \eqref{optimization_problem} are differentiable, algorithm (\ref{distributed project_v2_x}) becomes an ordinary differential equation and the strict convexity requirement of the objective functions can be relaxed, still with our proposed algorithm and technique.

\begin{corollary}\label{theorem_smooth_convergence}
Consider algorithm (\ref{distributed project_v2_x}). With 1) and 3) of Assumption \ref{Assumption}, if $f_i$ is differentiable  and convex on an open set containing $\Omega_i$,  and $\nabla f_i$ is Lipschitz continuous on $\Omega_i$ for $i\in\{1,\ldots,n\}$, then
\begin{itemize}
\item[($i$)] the solution $(x(t),\lambda(t),z(t))$ is bounded;
\item[($ii$)] the solution $(x(t),\lambda(t),z(t))$ is convergent and ${x}(t)$ converges to an optimal solution to problem \eqref{optimization_problem}.
\end{itemize}
\end{corollary}
\begin{IEEEproof}
($i$) The proof of ($i$) is similar to that of   Theorem \ref{theorem_project_convergence} ($i$). Hence, it is omitted.

($ii$) It follows from similar arguments in the proof of  Theorem \ref{theorem_project_convergence} ($ii$) that
\begin{eqnarray}
\frac{\mathrm d}{\mathrm{d}t}V (x,\lambda,z)\leq -\|\dot x\|^2 -(x-x^*)^{\rm T}(\nabla f({x})-\nabla f({x^*}))-\lambda^{\rm T}L\lambda\leq 0,
\end{eqnarray}
where $\dot x=P_{\Omega} [x  -  \nabla f({x})+\overline{W}^{\rm T} \lambda]- x$.

Let $\mathcal{R}=\big{\{} (x,\lambda,z)\in\Omega\times \mathbb{R}^{nm}\times \mathbb{R}^{nm}: 0=\frac{\mathrm d}{\mathrm{d}t}V(y,\lambda,z) \big{\}}\subset \big{\{} (x,\lambda,z)\in\Omega\times \mathbb{R}^{nm}\times \mathbb{R}^{nm}: L\lambda=0_{nm},(x-x^*)^{\rm T}(\nabla f({x})-\nabla f({x^*}))=0,\,\,0_{\sum_{i=1}^{n}q_i}=P_{\Omega} [x  -  \nabla f({x})+\overline{W}^{\rm T} \lambda]- x\big{\}}$.

Let $\mathcal M$ be the largest invariant subset of $\overline{\mathcal R}$.
It follows from the invariance principle (Theorem 2.41 of \cite{HC:2008}) that $(x(t),\lambda(t),z(t))\rightarrow\mathcal M$ as $t\rightarrow \infty$. Note that $\mathcal M$ is invariant. The trajectory $(\bar x(t), \bar \lambda(t),\bar z(t))\in \mathcal M$ for all $t\geq 0$ if $(\bar x(0), \bar \lambda(0),\bar z(0))=(\bar x_0, \bar \lambda_0,\bar z_0)\in\mathcal M$. Assume $(\bar x(t), \bar \lambda(t),\bar z(t))\in \mathcal M$ for all $t\geq 0$,
 $\dot {\bar x}(t)\equiv 0_{\sum_{i=1}^{q_i}}$ and $\dot {\bar z}(t)\equiv 0_{nm}$ and, hence, $\dot {\bar \lambda}(t)\equiv d-\overline {W }\bar x_0-L\bar z_0$.  Suppose $\dot {\bar \lambda}(t)\equiv d-\overline {W }\bar x_0-L\bar z_0\not=0_{nm}$, then ${\bar \lambda}(t)\rightarrow\infty$ as $t\rightarrow\infty$, which contradicts part ($i$). Hence, $\dot {\bar \lambda}(t)\equiv0_{nm}$ and $\mathcal M\subset \big{\{} (x,\lambda,z)\in\Omega\times \mathbb{R}^{nm}\times \mathbb{R}^{nm}: P_{\Omega} [x  -  \nabla f({x})+\overline{W}^{\rm T} \lambda]- x=0_{\sum_{i=1}^{n}q_i},\, d-\overline{W}x-Lz=0_{nm},\,L\lambda=0_{nm}\big{\}}$.

Take any $(\bar x,\bar {\lambda},\bar z)\in\mathcal M$.  Obviously, $(\bar x,\bar {\lambda},\bar z)$ is an equilibrium point of algorithm (\ref{distributed project_v2_x}). Define a new function $\bar V(x,\lambda,z)$ by replacing $( x^*, \lambda^*, z^*)$ with $(\bar x,\bar {\lambda},\bar z)$ in $V(x,\lambda,z)$. It follows from similar arguments in the proof of  Lemma \ref{differential_project_eqn} that $\frac{\mathrm d}{\mathrm {d}t}\bar V(x,\lambda,z)\leq 0$.
Hence, $(\bar x,\bar \lambda,\bar z)$ is Lyapunov stable. By Proposition 4.7 of \cite{HC:2008}, there exists $(\tilde x,\tilde \lambda,\tilde z)\in\mathcal M$ such that $ (x(t),\lambda (t),z(t))\rightarrow (\tilde x,\tilde \lambda,\tilde z)$ as $t\rightarrow\infty$.  Since $(\tilde x,\tilde \lambda,\tilde z)\in\mathcal M$ is an equilibrium point of algorithm (\ref{distributed project_v2_x}), $\tilde x$ is an optimal solution to problem \eqref{optimization_problem} by Theorem \ref{equivalent_project_conditions}.
\end{IEEEproof}

\section{Numerical Experiments}\label{sec_appl}
In this section, we give two numerical examples to illustrate the effectiveness of the proposed algorithms.

\subsection{Nonsmooth Optimization Problem}\label{sec_NS}
%We revise problem (\ref{QP_problem}) by adding nonsmooth terms to the objective functions,
Consider the following nonsmooth optimization problem
\begin{eqnarray}\label{nonsmooth_programming}
% \nonumber % Remove numbering (before each equation)
  \min  f(x) ,\quad f(x)= \sum_{i=1}^{n}(|x_i|^2+|x_i|),\quad \sum_{i=1}^{n}A_ix_i=d_0,\quad |x_i|\leq 1,\quad i\in\{1,\ldots,n\},
\end{eqnarray}
where $x=[x_1,\ldots,x_n]^{\rm T}\in\mathbb{R}^{n}$ and $x_i\in\mathbb{R}$. %Let  $n$, $A$, $d_0 $, $\omega_i$, $d_i $, $i\in\{1,\ldots,10\}$, be as defined in Section \ref{Quadratic_Programming}.
Let $n=10$, $d_0 = [3,2]^{\rm T}$,  $d_i = [0.3,0.2]^{\rm T}$, $i\in\{1,\ldots,10\}$, and
$$
A=\begin{bmatrix}
  1 & 1 & 1 & 0 & 0 & 1 & 1 & 1 & 0 & 0\\
  1 & 0 & 0 & 1 & 1 & 1 & 0 & 0 & 1 & 1
  \end{bmatrix}.
  $$

\begin{figure}
  \centering
  % Requires \usepackage{graphicx}
  \includegraphics[width=9 cm, height = 6 cm]{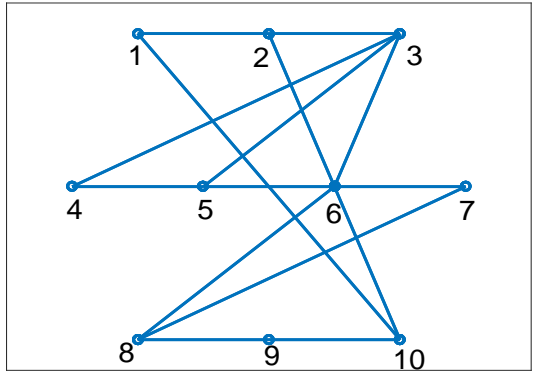}
  %\caption{Graph Topology of Optimization Algorithm in Section \ref{Quadratic_Programming}}
  \caption{Graph Topology}
  \label{fig_topo_10}
\end{figure}

\begin{figure}
  \centering
  % Requires \usepackage{graphicx}
  \includegraphics[width=12 cm, height = 8 cm]{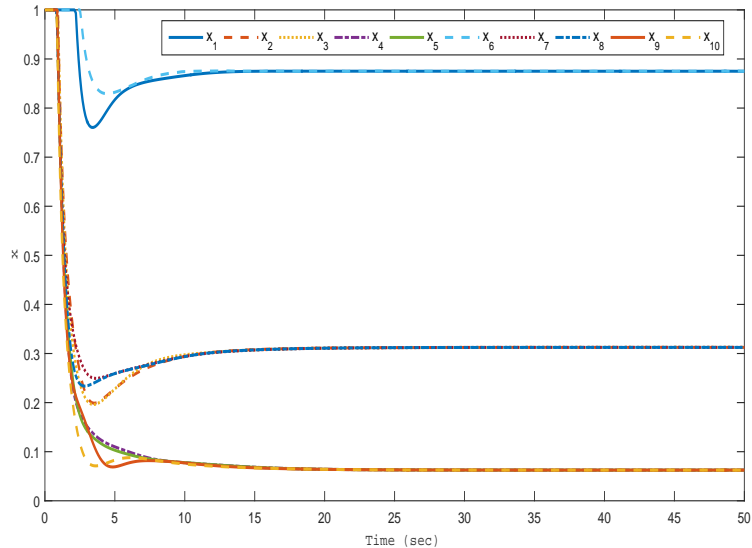}
  \caption{Trajectories of estimates for $x$ versus time  under algorithm \eqref{distributed feedback_v2_x} for problem \eqref{nonsmooth_programming}}
  \label{X_NS}
\end{figure}

\begin{figure}
  \centering
  % Requires \usepackage{graphicx}
  \includegraphics[width=12 cm, height = 8 cm]{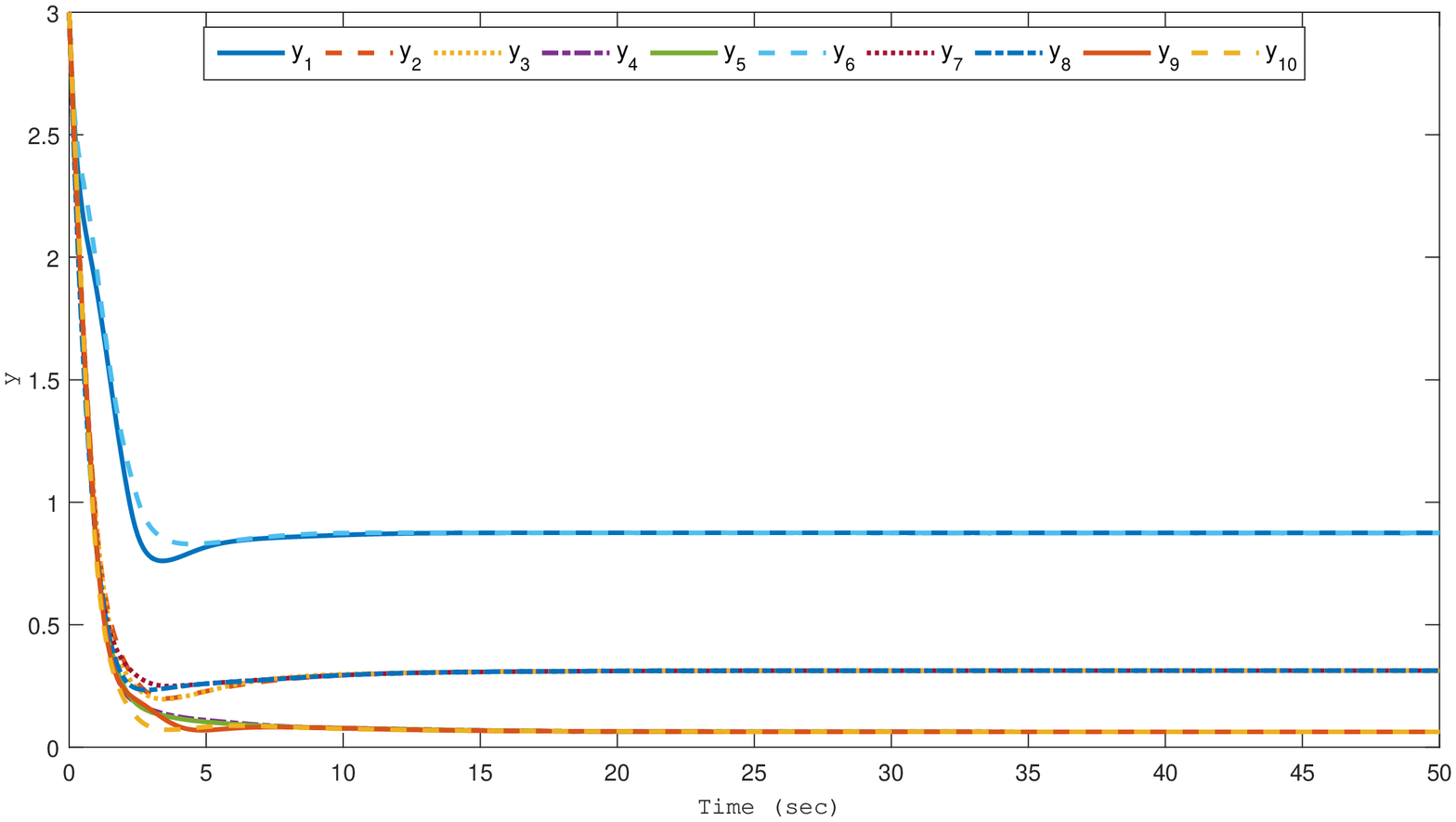}
  \caption{Trajectories of the auxiliary variable $y$ versus time  under algorithm \eqref{distributed feedback_v2_x} for problem \eqref{nonsmooth_programming}}
  \label{Y_NS}
\end{figure}

\begin{figure}
  \centering
  % Requires \usepackage{graphicx}
  \includegraphics[width=12 cm, height = 8 cm]{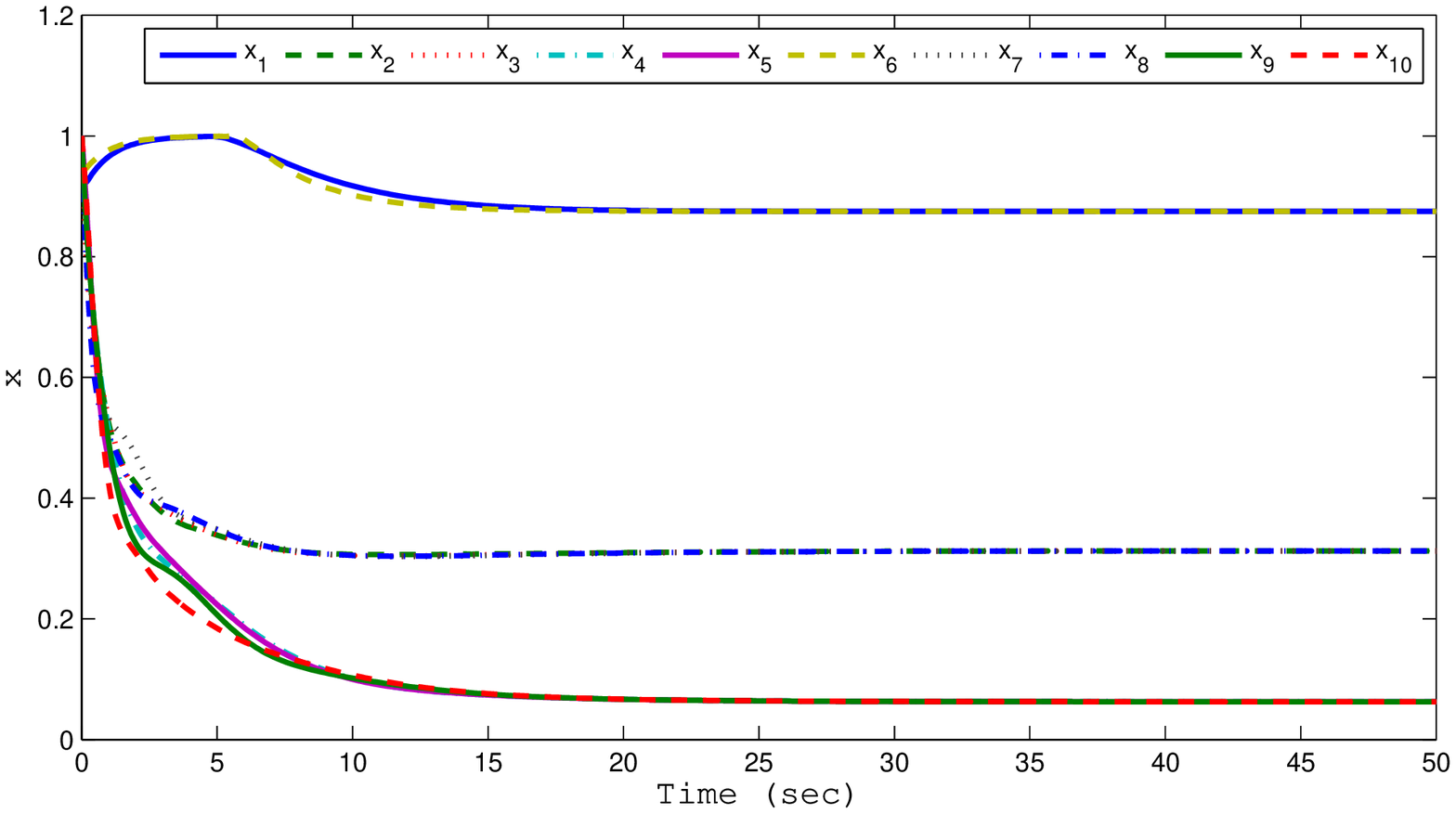}
  \caption{Trajectories of estimates for $x$ versus time  under algorithm \eqref{distributed project_v2_x} for problem \eqref{nonsmooth_programming}}
  \label{X_NS_PJ}
\end{figure}

\begin{figure}
  \centering
  % Requires \usepackage{graphicx}
  \includegraphics[width=12 cm, height = 8 cm]{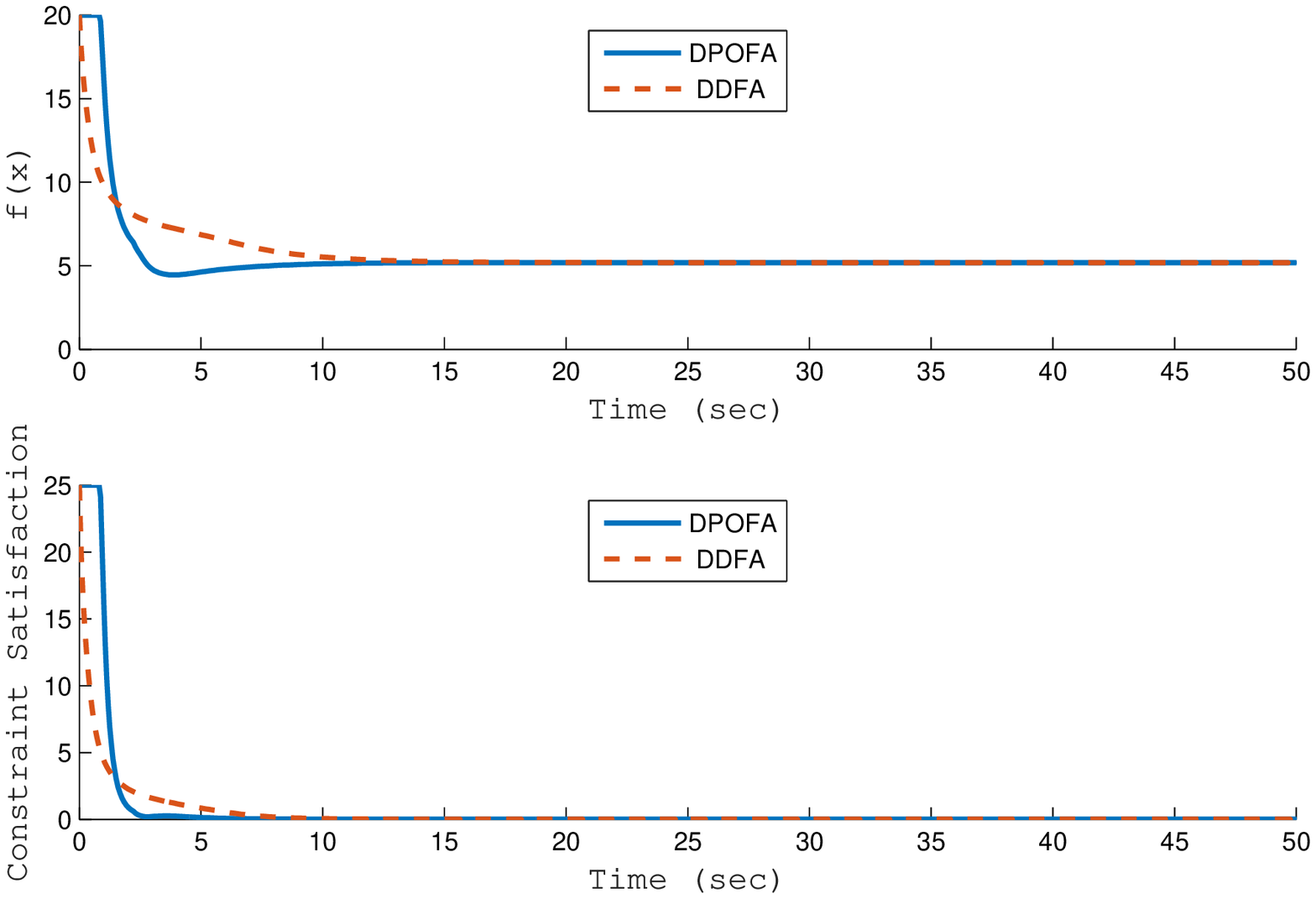}
  \caption{Objective functions $f(x)$ and constraints $\|Wx-d_0\|^2$ versus time  under algorithms \eqref{distributed feedback_v2_x} and \eqref{distributed project_v2_x} for problem \eqref{nonsmooth_programming}}
  \label{objective_constraint_NS}
\end{figure}

\begin{figure}
  \centering
  % Requires \usepackage{graphicx}
  \includegraphics[width=12 cm, height = 8 cm]{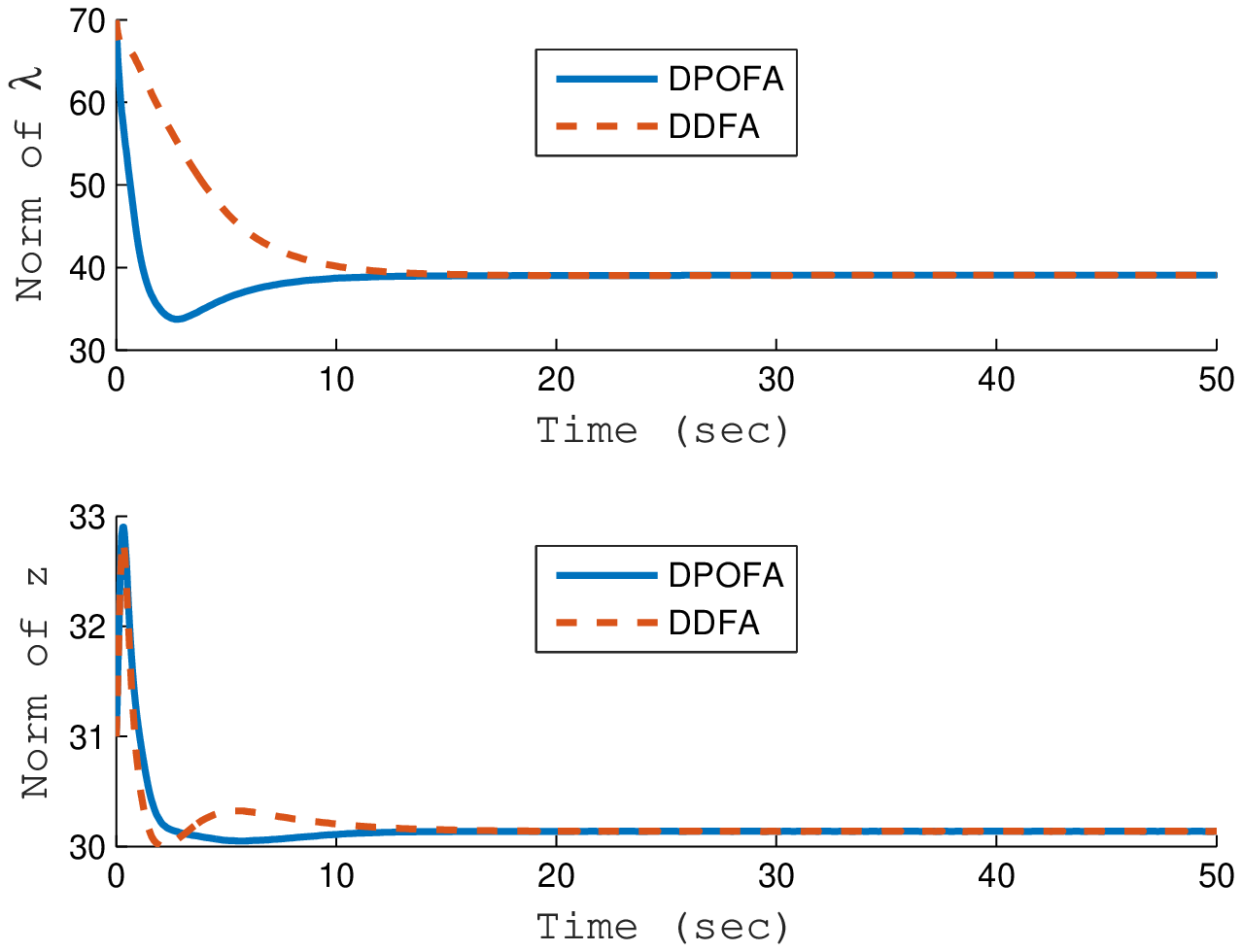}
  \caption{$\|\lambda\|^2$ and $\|z\|^2$ versus time  under algorithms \eqref{distributed feedback_v2_x} and \eqref{distributed project_v2_x} for problem \eqref{nonsmooth_programming}}
  \label{LBD_Z_NS}
\end{figure}

Take a network of 10 agents interacting over a graph $\mathcal G$ to solve this problem. The  information sharing graph  $\mathcal G$ of the optimization algorithms is given in Fig. \ref{fig_topo_10}. Some simulation results using  DPOFA algorithm \eqref{distributed feedback_v2_x} proposed in Section \ref{project_output} and DDFA algorithm \eqref{distributed project_v2_x} proposed in Section \ref{sec:alg_proj} are shown in Figs. \ref{X_NS}-\ref{LBD_Z_NS}.

Figs. \ref{X_NS} and \ref{Y_NS} show the trajectories of estimates for $x$ and the auxiliary variable $y$ versus time under DPOFA algorithm \eqref{distributed feedback_v2_x} proposed in Section \ref{project_output}, respectively, and
Fig. \ref{X_NS_PJ} depicts the trajectories of the estimates for $x$ versus time under DDFA algorithm \eqref{distributed project_v2_x} proposed in Section \ref{sec:alg_proj}.
Algorithm \eqref{distributed feedback_v2_x} proposed in Section \ref{project_output} uses an auxiliary variable $y$ and estimates  the optimal solution using $x=P_{\Omega}(y)$, while algorithm \eqref{distributed project_v2_x} proposed in Section \ref{sec:alg_proj} directly uses $x $ to estimate the solution. Both algorithms are able to find the optimal solution of the optimization problem. Figs. \ref{X_NS}-\ref{X_NS_PJ} indicate that the trajectories of $y$ may be out of the constraint set $\Omega$, but the  trajectories of $x$ stays in the constraint set $\Omega$.

Fig. \ref{objective_constraint_NS} gives the trajectories of the objective function $f(x)$ and constraint $\|Wx-d_0\|^2$ versus time under DPOFA algorithm \eqref{distributed feedback_v2_x} and  DDFA algorithm \eqref{distributed project_v2_x}, and demonstrates that the trajectories of $x$  converge to the equality constraint.
Furthermore, Fig. \ref{LBD_Z_NS} verifies the boundedness of the  trajectories of  auxiliary variables $\lambda$ and $z$.

In Fig. \ref{objective_constraint_NS}, the trajectories of $f(x)$ and $\|Wx-d_0\|^2$ versus time under DPOFA algorithm  show slow response speed at the beginning of the simulation.
This is because the change of $y$ in the algorithm may not generate the  changing behavior of $x=P_{\Omega}(y)$ when $y\notin\Omega$
%after the projection operation of $y$ to $\Omega$ at the beginning of simulation
(see the trajectories at the beginning of the simulation in Figs. \ref{X_NS} and \ref{Y_NS}).
Due to the indirect feedback effect on $x$ (changing  $x$ by controlling $y$) in DPOFA algorithm, the trajectory of variable $x$ may show slow changing behaviors in applications.

%\begin{figure}
%  \centering
%  % Requires \usepackage{graphicx}
%  \includegraphics[width=12 cm]{lbd1_ns.eps}
%  \caption{Trajectories of the auxiliary variable $\lambda_1$ versus time  under algorithm \eqref{distributed feedback_v2_x} for problem \eqref{nonsmooth_programming}}
%  \label{L1_NS}
%\end{figure}
%
%\begin{figure}
%  \centering
%  % Requires \usepackage{graphicx}
%  \includegraphics[width=12 cm]{lbd2_ns.eps}
%  \caption{Trajectories of the auxiliary variable $\lambda_2$ versus time  under algorithm \eqref{distributed feedback_v2_x} for problem \eqref{nonsmooth_programming}}
%  \label{L2_NS}
%\end{figure}
%
%
%\begin{figure}
%  \centering
%  % Requires \usepackage{graphicx}
%  \includegraphics[width=12 cm]{z1_ns.eps}
%  \caption{Trajectories of the auxiliary variable $z_1$ versus time  under algorithm \eqref{distributed feedback_v2_x} for problem \eqref{nonsmooth_programming}}
%  \label{Z1_NS}
%\end{figure}
%
%\begin{figure}
%  \centering
%  % Requires \usepackage{graphicx}
%  \includegraphics[width=12 cm]{z2_ns.eps}
%  \caption{Trajectories of the auxiliary variable $z_2$ versus time under algorithm \eqref{distributed feedback_v2_x} for problem \eqref{nonsmooth_programming}}
%  \label{Z2_NS}
%\end{figure}

%\begin{figure}
%  \centering
%  % Requires \usepackage{graphicx}
%  \includegraphics[width=12 cm, height = 8 cm]{f_x_NS_output.eps}
%  \caption{Objective functions $f(x)$ versus time  under algorithms \eqref{distributed feedback_v2_x} and \eqref{distributed project_v2_x} for problem \eqref{nonsmooth_programming}}
%  \label{objective_NS}
%\end{figure}

%\begin{figure}
%  \centering
%  % Requires \usepackage{graphicx}
%  \includegraphics[width=12 cm, height = 8 cm]{constraint_NS_control.eps}
%  \caption{$\|Wx-d\|^2$ versus time  under algorithms \eqref{distributed feedback_v2_x} and \eqref{distributed project_v2_x} for problem \eqref{nonsmooth_programming}}
%  \label{constraint_NS}
%\end{figure}

%\begin{figure}
%  \centering
%  % Requires \usepackage{graphicx}
%  \includegraphics[width=12 cm, height = 8 cm]{Z_NS.eps}
%  \caption{$\|z\|^2$ versus time  under algorithms \eqref{distributed feedback_v2_x} and \eqref{distributed project_v2_x} for problem \eqref{nonsmooth_programming}}
%  \label{Z_NS}
%\end{figure}

%\begin{figure}
%  \centering
%  % Requires \usepackage{graphicx}
%  \includegraphics[width=12 cm, height = 8 cm]{lambda_NS.eps}
%  \caption{$\|\lambda\|^2$ versus time  under algorithms \eqref{distributed feedback_v2_x} and \eqref{distributed project_v2_x} for problem \eqref{nonsmooth_programming}}
%  \label{lambda_NS}
%\end{figure}

\subsection{Multi-commodity Network Flow Problem}%{Min-Cost-Flow Problem}%Commodity Network Flow Problems
\label{Flow_Problem}
\subsubsection{Problem Description}
Consider a directed network consisting of a set $\mathcal{N}=\{1,\ldots,m\}$ of nodes and a set $\mathcal E$ of directed arcs. The flows on the arcs
are of $S$ different types (commodities). The flows of the $s$th type on the arc $(i,j)$ is denoted by $t_{i,j}(s)\in [l_{i,j}(s),u_{i,j}(s)]$, where $[l_{i,j}(s),u_{i,j}(s)]$ is the capacity constraint for $t_{i,j}(s)$. The flows must satisfy the conservation of flow and supply/demand constraints of the form
\begin{eqnarray}\label{demand_constraint}
% \nonumber % Remove numbering (before each equation)
   \sum_{j|(i,j)\in\mathcal E} t_{i,j}(s)-\sum_{j|(j,i)\in\mathcal E} t_{j,i}(s)= b_i(s),
\end{eqnarray}
where $ i\in\mathcal N,\, s\in\{1,\ldots,S\},$ $b_i(s)$ is the amount of flow of the $s$th type  entering the network at node $i$ ($b_i(s) > 0 $ indicates supply, and $b_i (s)< 0$ indicates demand). The supplies/demands $b_i(s)$ are given and satisfy $\sum_{i\in\mathcal N} b_i(s)=0$ for the  feasibility of the problem, which have been studied in the literature (see \cite{Bertsekas:JOTA:2008}).

The problem is to minimize $\sum_{(i,j)\in\mathcal E} f_{i,j}(t_{i,j}(1),\ldots,t_{i,j}(S))$ subject to the constraints \eqref{demand_constraint}, where $f_{i,j}:\mathbb R\rightarrow \mathbb R$ are continuous, strictly convex functions.

\subsubsection{Reformulation of Problem}
Let $n$ be the number of arcs in $\mathcal E$. We assign an index $k=k(i,j)\in\{1,\ldots,n\}$ to every arc $(i,j)\in\mathcal E$. We use $x_k=[t_{i,j}(1),\ldots, t_{i,j}(S)]^{\rm T}\in\mathbb R^S$ to denote the flow vector on arc $k$. Then constraint \eqref{demand_constraint} can be rewritten as
\begin{eqnarray}
% \nonumber % Remove numbering (before each equation)
  A\otimes I_S x =\sum_{k=1}^{n}A_k\otimes I_S x_k=b,
\end{eqnarray}
with $A$ as the vertex-edge incidence matrix of the graph, $A_k$ as the $k$th column of $A$, $x = [x_1^{\rm T},\ldots,x_n^{\rm T}]^{\rm T}$, $b_i=[b_i(1),\ldots,b_i(S)]^{\rm T}$, and $b = [b_1^{\rm T},\ldots,b_m^{\rm T}]^{\rm T}$.

Let $f^k(x_k)\triangleq f_{i,j}(t_{i,j}(1),\ldots, t_{i,j}(S))$, $\Omega_k=\prod_{s=1}^{S}  [l_{i,j}(s),u_{i,j}(s)]$, where $k$ is the index of arc $(i,j)\in\mathcal E$. The optimization problem can be reformulated as
\begin{eqnarray}\label{reformulate}
% \nonumber % Remove numbering (before each equation)
  \min f(x) = \sum_{k=1}^{n}f^k(x_k),\quad \sum_{k=1}^{n}A_k\otimes I_S x_k=b,\quad x_k\in \Omega_k,\quad k\in\{1,\ldots,n\}.
\end{eqnarray}

\subsubsection{Numerical Simulation}
\begin{figure}
  \centering
  % Requires \usepackage{graphicx}
  \includegraphics[width=12 cm, height = 8 cm]{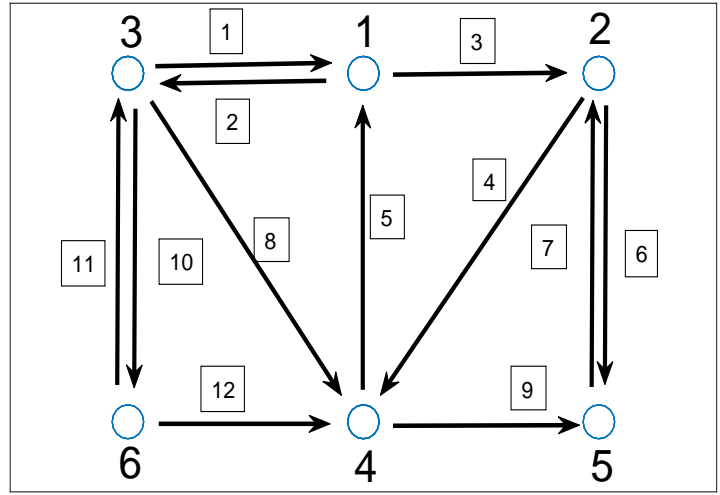}
  \caption{Graph Topology}
  \label{fig_topo1}
\end{figure}

Consider a network of 6 nodes and 12 arcs as shown in Fig. \ref{fig_topo1}, with $S=1$. Let $f^k(x_k) = \|x_k\|^2$, $x_k\in\mathbb R$ be the flows on arc $k$, and $\Omega_k=[0,10]$ for $k\in\{1,\ldots,12\}$. Problem \eqref{reformulate} can be formulated as
\begin{eqnarray}\label{prob_monotropic}
% \nonumber % Remove numbering (before each equation)
\min  \sum_{k=1}^{12} \|x_k\|^2, \quad Ax=\sum_{k=1}^{12}A_k x_k=b, \quad x_k\in \Omega_k,\quad k\in\{1,\ldots,12\},
\end{eqnarray}
where $b=[6,\,-7.2,\,-4.8,\,-9.6,\,8.4,\,7.2]^{\rm T}$ and $A_k$ is the $k$th column of the vertex-edge incidence matrix $A$ of the network in Fig. \ref{fig_topo1}.

\begin{figure}
  \centering
  % Requires \usepackage{graphicx}
  \includegraphics[width=12 cm, height = 8 cm]{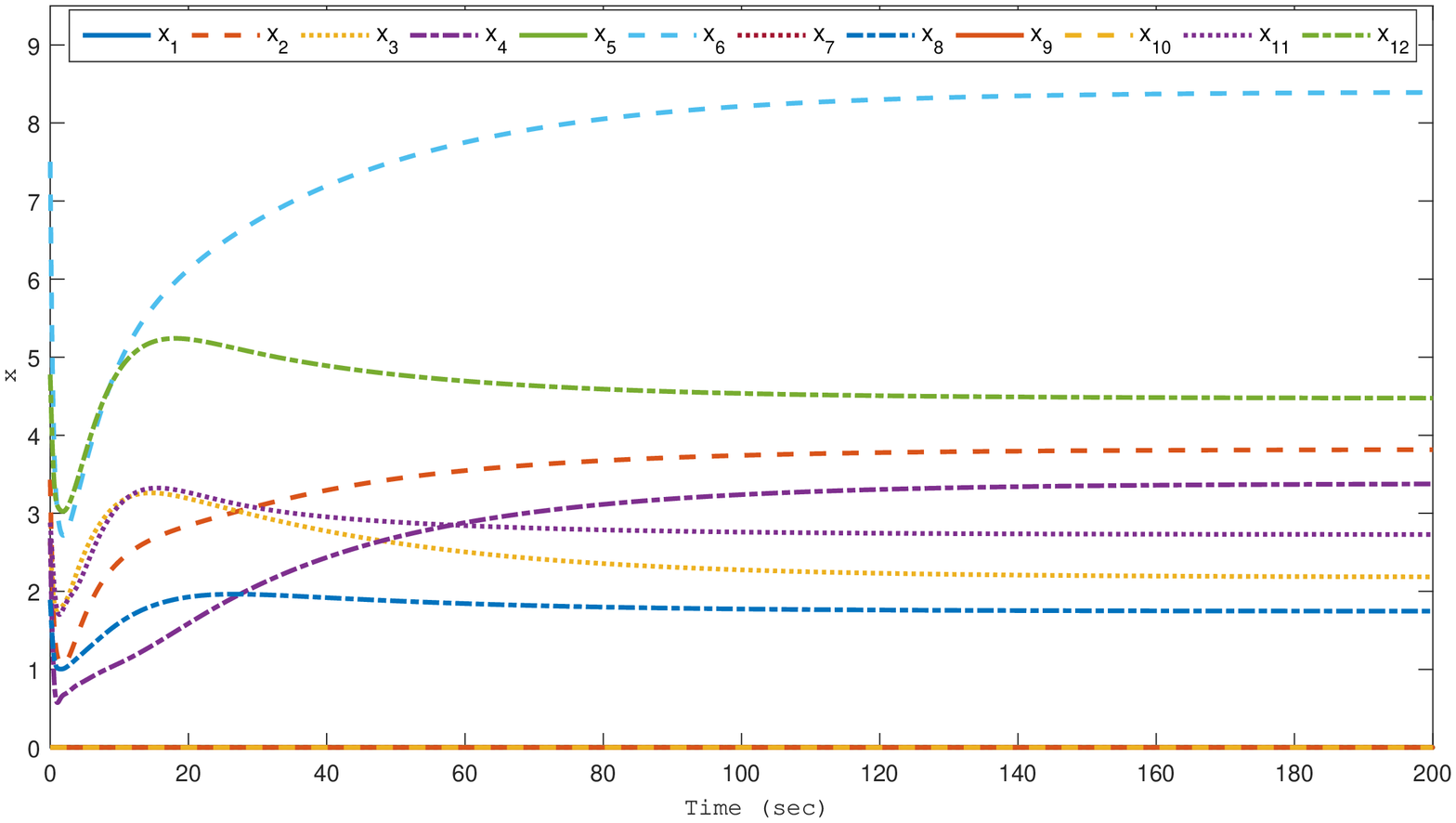}
  \caption{Trajectories of estimates for $x$ versus time  under algorithm \eqref{distributed feedback_v2_x} for problem \eqref{prob_monotropic}}
  \label{X_FP}
\end{figure}

\begin{figure}
  \centering
  % Requires \usepackage{graphicx}
  \includegraphics[width=12 cm, height = 8 cm]{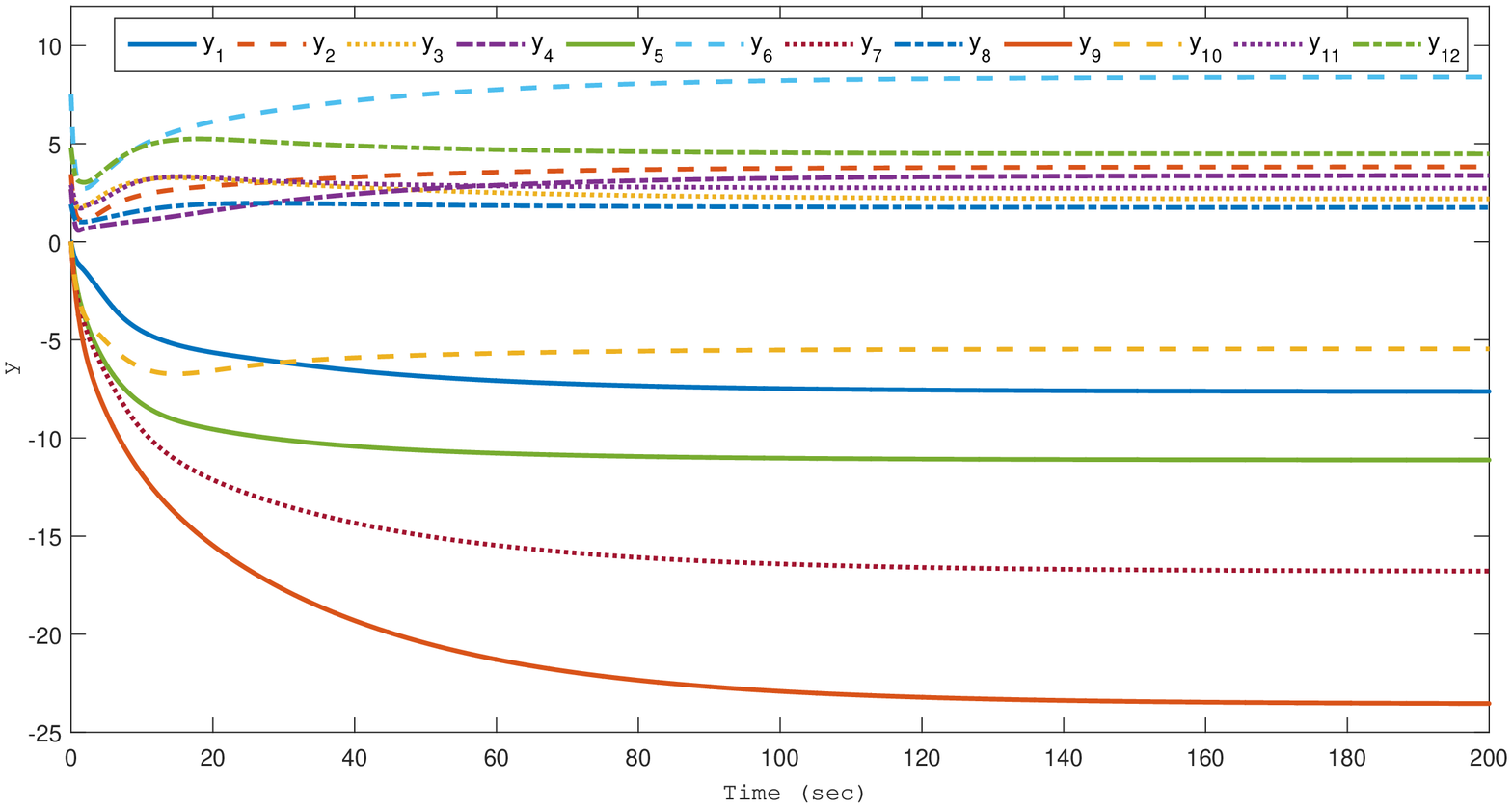}
  \caption{Trajectories of the auxiliary variable $y$ versus time under algorithm \eqref{distributed feedback_v2_x} for problem \eqref{prob_monotropic}}
  \label{Y_FP}
\end{figure}

\begin{figure}
  \centering
  % Requires \usepackage{graphicx}
  \includegraphics[width=12 cm, height = 8 cm]{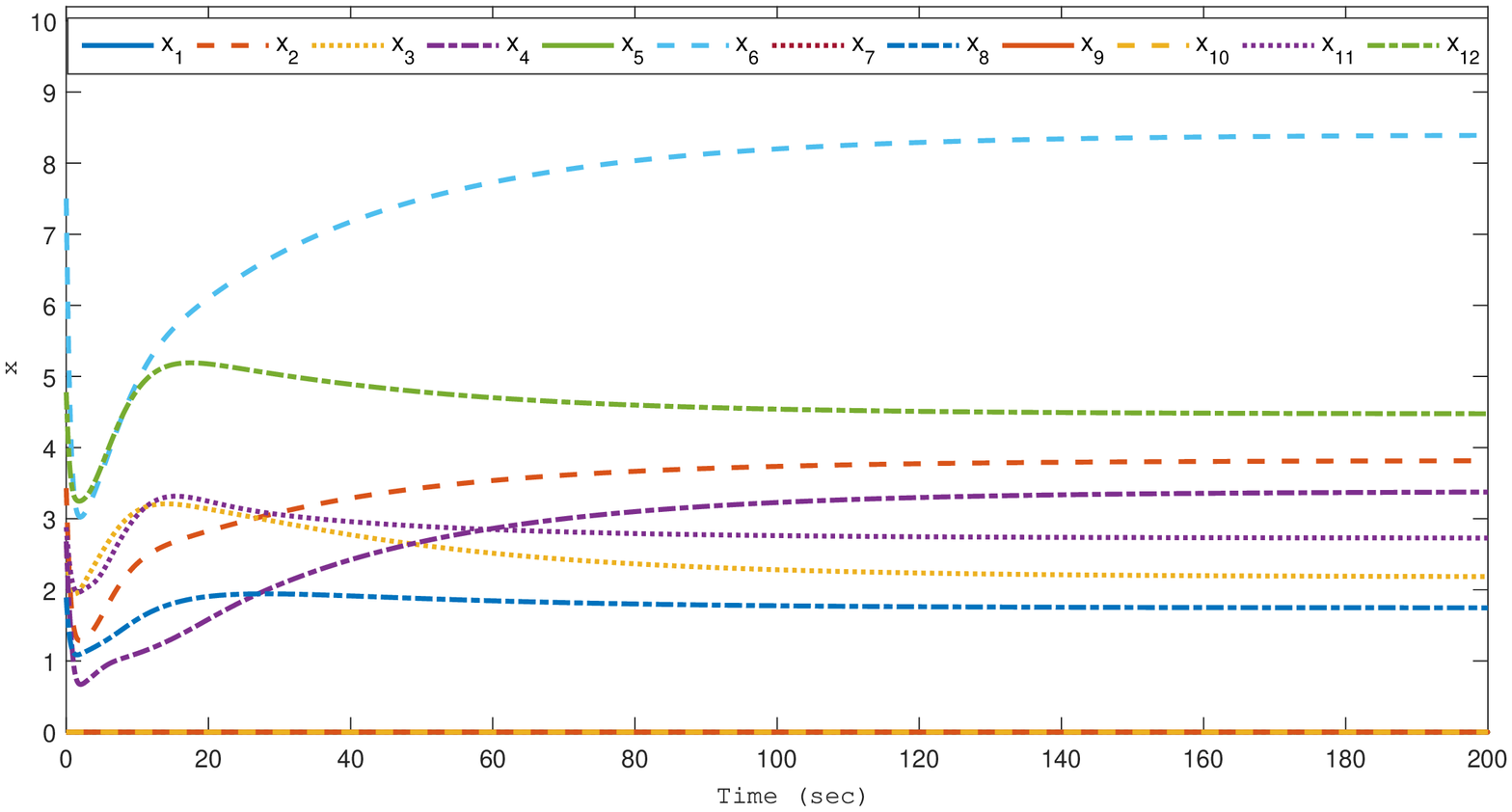}
  \caption{Trajectories of estimates for $x$ versus time  under algorithm \eqref{distributed project_v2_x} for problem \eqref{prob_monotropic}}
  \label{X_FP_PR}
\end{figure}

\begin{figure}
  \centering
  % Requires \usepackage{graphicx}
  \includegraphics[width=12 cm, height = 8 cm]{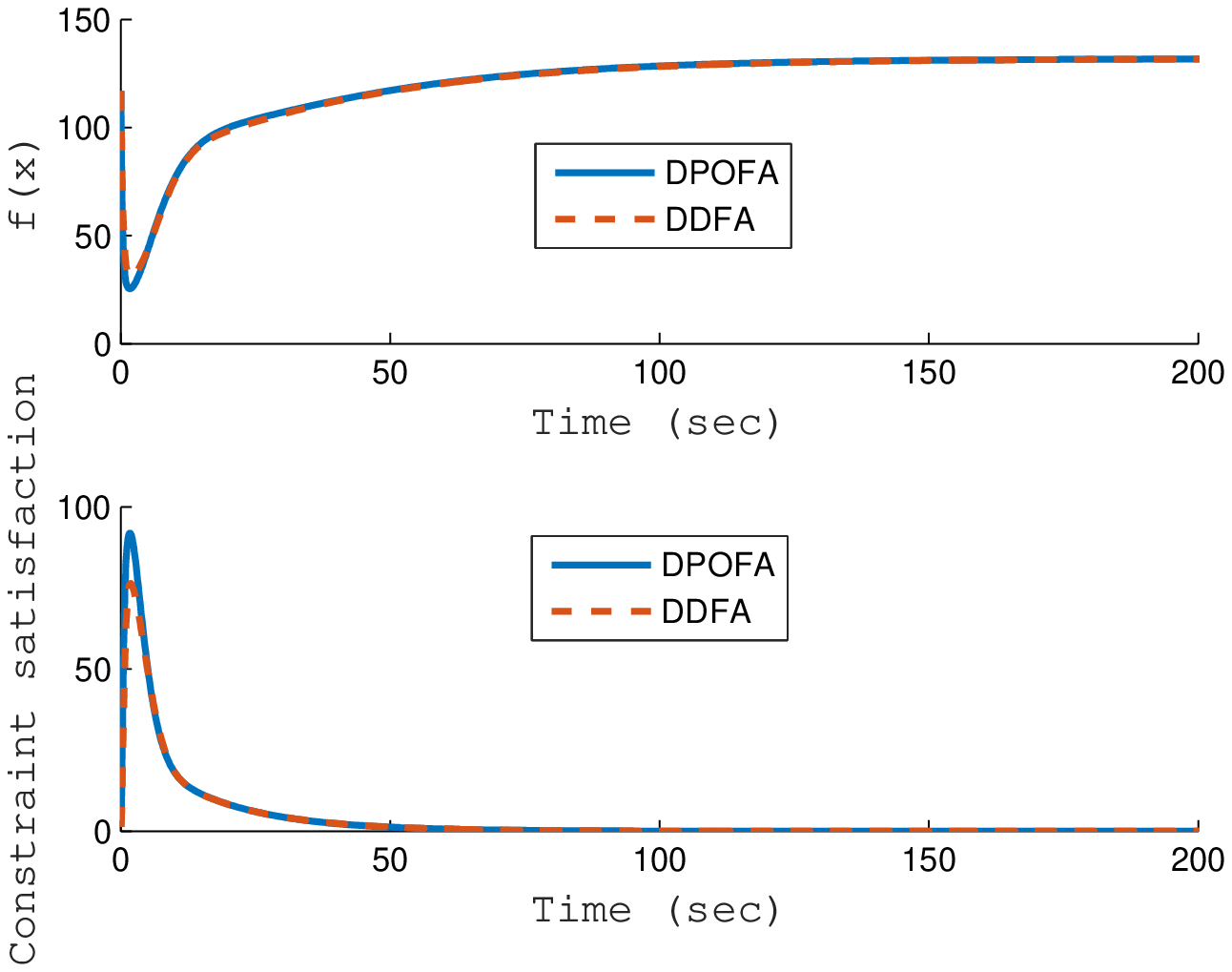}
  \caption{Objective functions $f(x)$ and constraints $\|Ax-b\|^2$ versus time  under algorithms \eqref{distributed feedback_v2_x} and \eqref{distributed project_v2_x} for problem \eqref{prob_monotropic}}
  \label{objective_constraint_FLOW}
\end{figure}

\begin{figure}
  \centering
  % Requires \usepackage{graphicx}
  \includegraphics[width=12 cm, height = 8 cm]{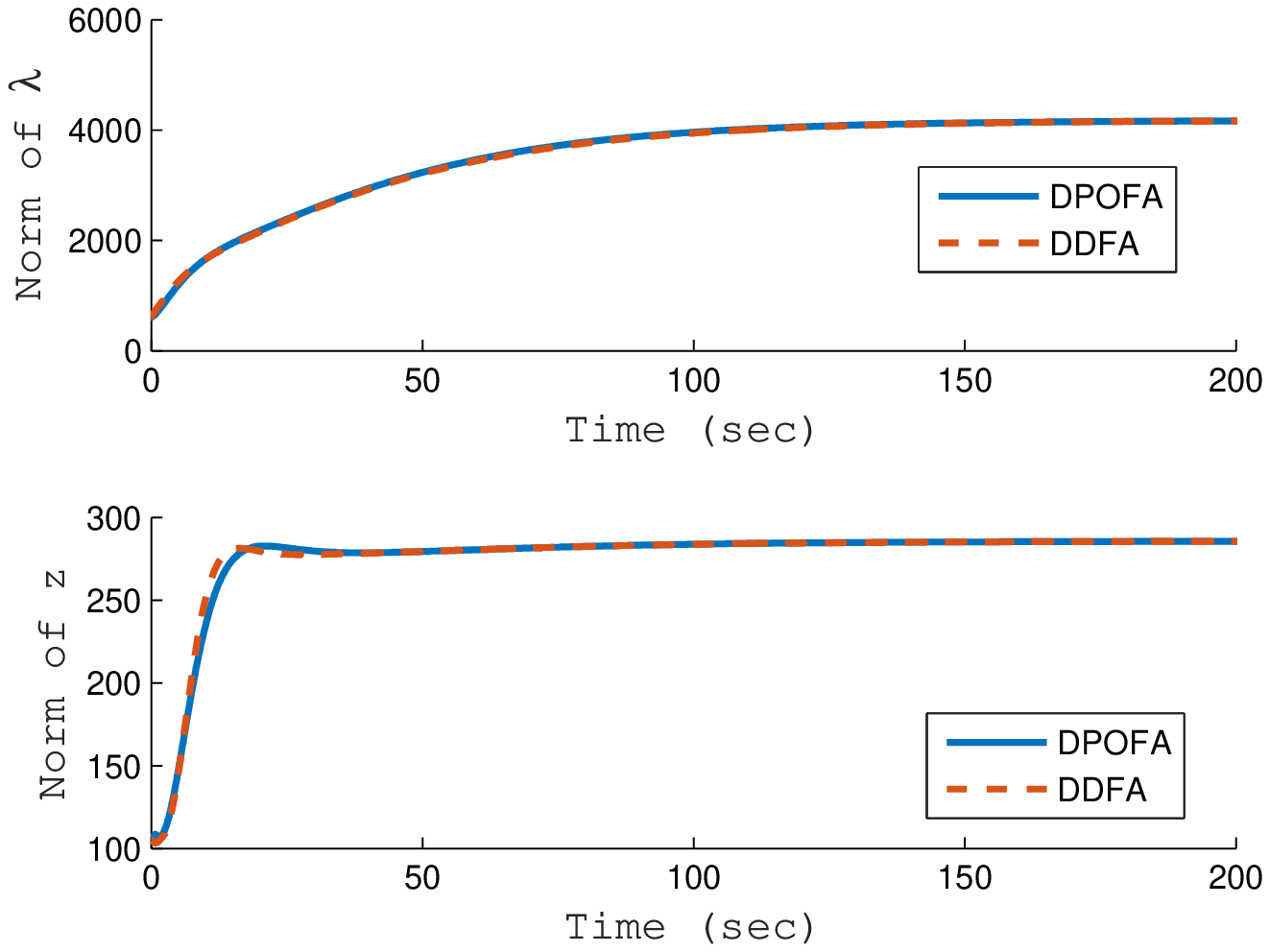}
  \caption{$\|\lambda\|^2$ and $\|z\|^2$ versus time  under algorithms \eqref{distributed feedback_v2_x} and \eqref{distributed project_v2_x} for problem \eqref{prob_monotropic}}
  \label{LBD_Z_FLOW}
\end{figure}

Simulation results using DPOFA algorithm \eqref{distributed feedback_v2_x} proposed in Section \ref{project_output} and DDFA algorithm \eqref{distributed project_v2_x} proposed in Section \ref{sec:alg_proj} are shown in Figs. \ref{X_FP}-\ref{LBD_Z_FLOW}.

Fig. \ref{X_FP} shows the trajectories of estimates for $x$ versus time of DPOFA algorithm \eqref{distributed feedback_v2_x} proposed in Section \ref{project_output}, while
Fig. \ref{Y_FP} shows those of the auxiliary variables $y$ versus time of DPOFA algorithm \eqref{distributed feedback_v2_x} proposed in Section \ref{project_output}.
Fig. \ref{X_FP_PR} exhibits the trajectories of estimates for $x$ versus time of DDFA algorithm \eqref{distributed project_v2_x} proposed in Section \ref{sec:alg_proj}.
Algorithm \eqref{distributed feedback_v2_x} proposed in Section \ref{project_output} uses an auxiliary variable $y$ and estimates  the optimal solution using $x=P_{\Omega}(y)$, while algorithm \eqref{distributed project_v2_x} proposed in Section \ref{sec:alg_proj} directly uses $x $ to estimate the solution. Both algorithms are able to find the optimal solution of the optimization problem. Figs. \ref{X_FP}-\ref{X_FP_PR} indicate that the trajectories of $y$ may be out of the constraint set $\Omega$, but the  trajectories of $x$  stay in the constraint set $\Omega$.

Fig. \ref{objective_constraint_FLOW} depicts the trajectories of the objective function $f(x)$ and constraint $\|Ax-b\|^2$ versus time under DPOFA algorithm \eqref{distributed feedback_v2_x} and DDFA algorithm \eqref{distributed project_v2_x}, where the trajectories of $x$  converge to the equality constraint. Fig. \ref{LBD_Z_FLOW}  illustrates the boundedness of the  trajectories of  auxiliary variables $\lambda$ and $z$.  Clearly, both algorithms can solve the problem.% and the convergence rates are quite similar.

%\begin{eqnarray}
%A=\begin{bmatrix}
%-1 & 1 &  1 & 0 &  -1 & 0 & 0  & 0  & 0         & 0 &  0 &  0\\
%0  & 0 & -1 & 1  & 0  & -1 & 1 &  0  & 0       & 0 & 0 &  0\\
%1 &  -1 & 0 & 0  & 0 &  0 &  0  & 1 &  0       &  1 &  -1 & 0\\
%0 &  0 &  0 & -1 & 1  & 0 &  0 &   -1 & 1       &  0 &  0 &  -1\\
%0 &  0  & 0 &  0 &  0 &  1 & -1 & 0  & -1      & 0 &  0 &  0\\
%0 &  0 & 0 &  0 &  0  & 0 & 0 &  0 &  0        &  -1 & 1  & 1
%\end{bmatrix}.
%\end{eqnarray}

%\begin{figure}
%  \centering
%  % Requires \usepackage{graphicx}
%  \includegraphics[width=12 cm, height = 8 cm]{flow_objective.eps}
%  \caption{Objective functions $f(x)$ versus time under algorithms \eqref{distributed feedback_v2_x} and \eqref{distributed project_v2_x} for problem \eqref{prob_monotropic}}
%  \label{objective_FP}
%\end{figure}
%
%\begin{figure}
%  \centering
%  % Requires \usepackage{graphicx}
%  \includegraphics[width=12 cm, height = 8 cm]{flow_constraint.eps}
%  \caption{$\|Wx-d\|^2$ versus time under algorithms \eqref{distributed feedback_v2_x} and \eqref{distributed project_v2_x} for problem \eqref{prob_monotropic}}
%  \label{constriant_FP}
%\end{figure}
%
%\begin{figure}
%  \centering
%  % Requires \usepackage{graphicx}
%  \includegraphics[width=12 cm, height = 8 cm]{z_flow.eps}
%  \caption{$\|z\|^2$ versus time  under algorithms \eqref{distributed feedback_v2_x} and \eqref{distributed project_v2_x} for problem \eqref{prob_monotropic}}
%  \label{Z_flow}
%\end{figure}
%...

\section{Conclusions}\label{conclusion}

In this paper, the distributed design for the extended monotropic optimization (EMO) problem has been addressed, which is related to various applications in large-scale optimization and evolutionary computation.  In this paper, two novel distributed continuous-time algorithms using  projected output feedback design and derivative feedback design have been proposed to solve this problem in multi-agent networks.
The design of the algorithms  has used the decomposition of problem constraints and distributed techniques. Based on  stability theory and invariance principle for differential inclusions, the convergence properties of the proposed algorithms have been established. The trajectories of all the agents have been proved to be bounded and convergent  to the optimal solution with any initial condition in mathematical and numerical ways.

The distributed EMO problem definitely deserves more efforts because of its broad range of applications. In  the future, distributed EMO problems with more complicated situations   such as  parameter uncertainties and
online concerns  will be further investigated.

%--------------------------------------------------------
%This is an alternative way to create a reference list: Put all %the cited references below into Reference.bib
%You can create such a bib file by yourself
%
\bibliographystyle{IEEEtran}
\bibliography{Reference}

% Generated by IEEEtran.bst, version: 1.14 (2015/08/26)
\begin{thebibliography}{10}
\providecommand{\url}[1]{#1}
\csname url@samestyle\endcsname
\providecommand{\newblock}{\relax}
\providecommand{\bibinfo}[2]{#2}
\providecommand{\BIBentrySTDinterwordspacing}{\spaceskip=0pt\relax}
\providecommand{\BIBentryALTinterwordstretchfactor}{4}
\providecommand{\BIBentryALTinterwordspacing}{\spaceskip=\fontdimen2\font plus
\BIBentryALTinterwordstretchfactor\fontdimen3\font minus
  \fontdimen4\font\relax}
\providecommand{\BIBforeignlanguage}[2]{{%
\expandafter\ifx\csname l@#1\endcsname\relax
\typeout{** WARNING: IEEEtran.bst: No hyphenation pattern has been}%
\typeout{** loaded for the language `#1'. Using the pattern for}%
\typeout{** the default language instead.}%
\else
\language=\csname l@#1\endcsname
\fi
#2}}
\providecommand{\BIBdecl}{\relax}
\BIBdecl

\bibitem{LW:TAC:2015}
Q.~Liu and J.~Wang, ``A second-order multi-agent network for bound-constrained
  distributed optimization,'' \emph{IEEE Transaction on Automatic Control},
  vol.~60, no.~12, pp. 3310--3315, 2015.

\bibitem{BX:2009}
W.~Bian and X.~Xue, ``Subgradient-based neural networks for nonsmooth nonconvex
  optimization problems,'' \emph{IEEE Transactions on Neural Networks},
  vol.~20, no.~6, pp. 1024--1038, 2009.

\bibitem{YHL:SCL:2015}
P.~Yi, Y.~Hong, and F.~Liu, ``Distributed gradient algorithm for constrained
  optimization with application to load sharing in power systems,''
  \emph{Systems \& Control Letters}, vol.~83, pp. 45--52, 2015.

\bibitem{YHX:TNNLS:2016}
D.~Yuan, D.~W.~C. Ho, and S.~Xu, ``Zeroth-order method for distributed
  optimization with approximate projections,'' \emph{IEEE Transactions on
  Neural Networks and Learning Systems}, vol.~27, no.~2, pp. 284--294, 2016.

\bibitem{NOP:2010}
A.~Nedic, A.~Ozdaglar, and P.~A. Parrilo, ``Constrained consensus and
  optimization in multi-agent networks,'' \emph{IEEE Transactions on Automatic
  Control}, vol.~55, pp. 922--938, 2010.

\bibitem{GC:2014}
B.~Gharesifard and J.~Cort\'{e}s, ``Distributed continuous-time convex
  optimization on weight-balanced digraphs,'' \emph{IEEE Transactions on
  Automatic Control}, vol.~59, pp. 781--786, 2014.

\bibitem{SJH:TAC:2013}
G.~Shi, K.~Johansson, and Y.~Hong, ``Reaching an optimal consensus: Dynamical
  systems that compute intersections of convex sets,'' \emph{IEEE Transactions
  on Automatic Control}, vol.~55, no.~3, pp. 610--622, 2013.

\bibitem{LYW:NNLS:2016}
Q.~Liu, S.~Yang, and J.~Wang, ``A collective neurodynamic approach to
  distributed constrained optimization,'' \emph{IEEE Transactions on Neural
  Networks and Learning Systems}, to appear.

\bibitem{TH:1986}
D.~Tank and J.~Hopfield, ``Simple neural optimization networks: An a/d
  converter, signal decision circuit, and a linear programming circuit,''
  \emph{IEEE Transactions on Circuits and Systems}, vol.~33, no.~5, pp.
  533--541, 2013.

\bibitem{CHZ:1999}
E.~Chong, S.~Hui, and S.~Zak, ``An analysis of a class of neural networks for
  solving linear programming problems,'' \emph{IEEE Transactions on Automatic
  Control}, vol.~28, no.~3, pp. 36--73, 2008.

\bibitem{LW:MAY2013}
Q.~Liu and J.~Wang, ``A one-layer projection neural network for nonsmooth
  optimization subject to linear equalities and bound constraints,'' \emph{IEEE
  Transactions on Neural Networks and Learning Systems}, vol.~24, no.~5, pp.
  812--824, 2013.

\bibitem{PSL:2003}
J.~Park, K.~Lee, J.~Shin, and K.~Y. Lee, ``Economic load dispatch for
  non-smooth cost functions using particle swarm optimization,'' in \emph{IEEE
  Power Engineering Society 2003 General Meeting}, Toronto, Canada, 2003, pp.
  938--943.

\bibitem{RM:2015}
E.~Ram\'{i}rez-Llanos and S.~Mart\'{i}nez, ``Distributed and robust resource
  allocation algorithms for multi-agent systems via discrete-time iterations,''
  in \emph{Proc. IEEE Conf. Decision Control}, Osaka, Japan, 2015, pp.
  1390--1396.

\bibitem{arXiv:YHL}
P.~Yi, Y.~Hong, and F.~Liu, ``Initialization-free distributed algorithms for
  optimal resource allocation with feasibility constraints and its application
  to economic dispatch of power systems,'' \emph{{arXiv:1510.08579}}, 2015.

\bibitem{FNQ:2004}
P.~N. M.~Forti and M.~Quincampoix, ``Generalized neural network for nonsmooth
  nonlinear programming problems,'' \emph{IEEE Transactions on Circuits and
  Systems-I}, vol.~51, no.~9, pp. 1741--1754, 2004.

\bibitem{Rockafellar:1984}
R.~Rockafellar, \emph{Network Flows and Monotropic Optimization}.\hskip 1em
  plus 0.5em minus 0.4em\relax New York: Wiley, 1984.

\bibitem{Rockafellar1985}
------, ``Monotropic programming: A generalization of linear programming and
  network programming. in: Convexity and duality in optimization,''
  \emph{Lecture Notes in Economics and Mathematical Systems}, vol. 256, pp.
  226--237, 1985.

\bibitem{Bertsekas:JOTA:2008}
D.~P. Bertsekas, ``Extended monotropic programming and duality,'' \emph{Journal
  of Optimization Theory and Applications}, vol. 139, no.~2, pp. 209--225,
  2008.

\bibitem{CDZ:MPA:2015}
N.~Chatzipanagiotis, D.~Dentcheva, and M.~M. Zavlanos, ``An augmented
  lagrangian method for distributed optimization,'' \emph{Mathematical
  Programming Series A}, vol. 152, no.~1, pp. 405--434, 2015.

\bibitem{GR2001}
C.~Godsil and G.~F. Royle, \emph{Algebraic Graph Theory}.\hskip 1em plus 0.5em
  minus 0.4em\relax New York: Springer-Verlag, 2001.

\bibitem{AC:1984}
J.~P. Aubin and A.~Cellina, \emph{Differential Inclusions}.\hskip 1em plus
  0.5em minus 0.4em\relax Berlin, Germany: Springer-Verlag, 1984.

\bibitem{Clarke:1983}
F.~H. Clarke, \emph{Optimization and Nonsmooth Analysis}.\hskip 1em plus 0.5em
  minus 0.4em\relax New York: Wiley, 1983.

\bibitem{BC:1999}
A.~Bacciotti and F.~Ceragioli, ``Stability and stabilization of discontinuous
  systems and nonsmooth {L}yapunov functions,'' \emph{ESAIM: Control,
  Optimisation and Calculus of Variations}, vol.~4, pp. 361--376, 1999.

\bibitem{Cortez:2008}
J.~Cort\'{e}s, ``Discontinuous dynamical systems,'' \emph{IEEE Control Systems
  Magazine}, vol.~44, no.~11, pp. 1995--2006, 1999.

\bibitem{ZSN:2003}
Z.~Denkowski, S.~Mig\'{o}rski, and N.~S. Papageorgiou, \emph{An Introduction to
  Nonlinear Analysis: Theory}.\hskip 1em plus 0.5em minus 0.4em\relax New York,
  NY: Springer-Verlag New York Inc., 2003.

\bibitem{KS:1982}
D.~Kinderlehrer and G.~Stampacchia, \emph{An Introduction to Variational
  Inequalities and Their Applications}.\hskip 1em plus 0.5em minus 0.4em\relax
  New York: Academic, 1982.

\bibitem{QLX:A:2016}
Z.~Qiu, S.~Liu, and L.~Xie, ``Distributed constrained optimal consensus of
  multi-agent systems,'' \emph{Automatica}, vol.~68, pp. 209--215, 2016.

\bibitem{arXiv:0905.4745}
M.~Tygert, ``A fast algorithm for computing minimal-norm solutions to
  underdetermined systems of linear equations,'' \emph{{arXiv:0905.4745}},
  2009.

\bibitem{MLM:TAC:2015}
S.~Mou, J.~Liu, and A.~S. Morse, ``A distributed algorithm for solving a linear
  algebraic equation,'' \emph{IEEE Transactions on Automatic Control}, vol.~60,
  no.~11, pp. 2863--2878, 2015.

\bibitem{Ruszczynski:2006}
A.~Ruszczynski, \emph{Nonlinear Optimization}.\hskip 1em plus 0.5em minus
  0.4em\relax Princeton, New Jersey: Princeton University Press, 2006.

\bibitem{CC:network:2015}
A.~Cherukuri and J.~Cort\'{e}s, ``Distributed generator coordination for
  initialization and anytime optimization in economic dispatch,'' \emph{IEEE
  Transactions on Control of Network Systems}, vol.~2, no.~3, pp. 226--237,
  2015.

\bibitem{Strang:1993}
G.~Strang, ``The fundamental theorem of linear algebra,'' \emph{American
  Mathematical Monthly}, vol. 100, no.~9, pp. 848--855, 1993.

\bibitem{HC:2008}
W.~M. Haddad and V.~Chellaboina, \emph{Nonlinear Dynamical Systems and Control:
  A {L}yapunov-Based Approach}.\hskip 1em plus 0.5em minus 0.4em\relax
  Princeton, NJ: Princeton Univ. Press, 2008.

\end{thebibliography}

\end{document}